\documentclass[11pt]{article}%
\usepackage{amssymb}
\usepackage[all]{xy}
\usepackage[latin1]{inputenc}
\usepackage{amsfonts}
\usepackage{amsmath}
\usepackage{amssymb}
\usepackage{url}
\usepackage{hyperref}
\usepackage{graphicx}%
\providecommand{\U}[1]{\protect\rule{.1in}{.1in}}
\numberwithin{equation}{section}
\providecommand{\U}[1]{\protect\rule{.1in}{.1in}}
\providecommand{\U}[1]{\protect\rule{.1in}{.1in}}
\newtheorem{theo}{Theorem}[section]
\newtheorem{prop}[theo]{Proposition}

\newtheorem{defi}[theo]{Definition}

\begin{document}

\title{Fluctuations of Interacting Markov Chain Monte Carlo Methods}
\author{Bernard Bercu\thanks{Centre INRIA Bordeaux Sud-Ouest \& Institut de
Mathématiques de Bordeaux , Université Bordeaux, 351 cours de la Libération
33405 Talence cedex, France, Bernard.Bercu@math.u-bordeaux1.fr} , Pierre Del
Moral\thanks{Centre INRIA Bordeaux Sud-Ouest \& Institut de Mathématiques de
Bordeaux , Université Bordeaux, 351 cours de la Libération 33405 Talence
cedex, France, Pierre.Del-Moral@inria.fr} , Arnaud Doucet\thanks{Department of
Statistics, University of Oxford, 1 South Parks Road, Oxford OX1 3TG, United
Kingdom, Tel: +44(0)1865 282851 - Fax: +44(0)1865 272595,
Doucet@stats.ox.ac.uk}}
\date{1$^{\text{st}}$ version 30th Jan. 2008, 2$^{\text{nd}}$ revision 30th Dec.
2011. }
\maketitle

\begin{abstract}
We present a multivariate central limit theorem for a general class of
interacting Markov chain Monte Carlo algorithms used to solve nonlinear
measure-valued equations. These algorithms generate stochastic processes which
belong to the class of nonlinear Markov chains interacting with their
empirical occupation measures. We develop an original theoretical analysis
based on resolvent operators and semigroup techniques to analyze the
fluctuations of their occupation measures around their limiting values.


\emph{Mathematics Subject Classification} : \newline Primary: 47H20, 60G35,
60J85, 62G09; Secondary: 47D08, 47G10, 62L20.

\end{abstract}

\section{Introduction}

\subsection{Nonlinear measure-valued equations}

\label{nleqsec}

Let $(S^{(l)},\mathcal{S}^{(l)})_{l\geq0}$ be a sequence of measurable spaces.
For any $l\geq0$ we denote by $\mathcal{P}(S^{(l)})$ the set of probability
measures on $S^{(l)}$. Suppose we have a sequence of probability measures
$\pi^{(l)}\in\mathcal{P}(S^{(l)})$ such that, for any $l\geq1$, $\pi^{(l)}$
satisfies the following nonlinear measure-valued equation
\begin{equation}
\Phi^{\left(  l\right)  }(\pi^{(l-1)})=\pi^{(l)}\label{phi}%
\end{equation}
for some mappings $\Phi^{\left(  l\right)  }:\mathcal{P}(S^{(l-1)}%
)\rightarrow\mathcal{P}(S^{(l)})$. We will also use the convention
$\Phi^{\left(  0\right)  }(\pi^{(-1)})=\pi^{(0)}$.

In numerous scenarios, the probability measures $(\pi^{(l)})_{l\geq0}$ need to
be approximated numerically. Interacting particle methods have been previously
proposed to approximate these probability distributions \cite{arnaud,fk}.
However they suffer from several limitations detailed in \cite{brockwell},
\cite{dd08}. To bypass some of these limitations, an alternative class of
algorithms known as interacting Markov chain Monte Carlo (i-MCMC) methods has
been recently introduced in \cite{brockwell}, \cite{dd08}. The main objective
of this article is to present a multivariate Central Limit Theorem (CLT) for
these methods. This extends significantly our previous result established in
\cite{bdd08} which only applies to a restricted class of i-MCMC\ algorithms.

Before describing two general class of models where i-MCMC methods can be
used, we introduce the notation adopted in this paper.

\subsection{Notation and conventions}

\label{notat} We denote respectively by $\mathcal{M}(E)$, $\mathcal{M}_{0}(E)
$, $\mathcal{P}(E)$, and $\mathcal{B}(E)$, the set of all finite signed
measures on some measurable space $(E,\mathcal{E})$ equipped with some
$\sigma$-field $\mathcal{E }$, the convex subset of measures with null mass,
the subset of all probability measures, and finally the Banach space of all
bounded and measurable functions $f$ on $E$ equipped with the uniform norm
$\Vert f\Vert=\sup_{x\in E}{|f(x)|}$ and the Borel $\sigma$-field associated
to the supremum norm. We also denote by $\mathcal{B}_{1}(E)\subset
\mathcal{B}(E)$ the unit ball of functions $f\in\mathcal{B}(E)$ with $\Vert
f\Vert\leq1$, and by $\mbox{Osc}_{1}(E)$, the convex set of $\mathcal{E}%
$-measurable functions $f$ with oscillations less than one; that is,
\[
\mbox{osc}(f)=\sup{\{|f(x)-f(y)|\;;\;x,y\in E\}}\leq1
\]

We let $\mu(f)=\int~\mu(dx)~f(x)$ be the Lebesgue integral of a function
$f\in\mathcal{B}(E)$ with respect to a measure $\mu\in\mathcal{M}(E)$. We
slightly abuse the notation, and sometimes we denote by $\mu(A)=\mu(1_{A})$
the measure of a measurable subset $A\in\mathcal{E}$.

We recall that a bounded integral operator $M$ from a measurable space
$(E,\mathcal{E})$ into an auxiliary measurable space $(F,\mathcal{F})$ is an
operator $f\mapsto M(f)$ from $\mathcal{B}(F)$ into $\mathcal{B}(E)$ so that
the functions
\[
M(f)(x)=\int_{F}~M(x,dy)~f(y)~\in\mathbb{R}%
\]
are $\mathcal{E}$-measurable and bounded for any $f\in\mathcal{B}(F)$. By
Fubini's theorem, we recall that a bounded integral operator $M$ from a
measurable space $(E,\mathcal{E})$ into an auxiliary measurable space
$(F,\mathcal{F})$ also generates a dual operator $\mu\mapsto\mu M$ from
$\mathcal{M}(E)$ into $\mathcal{M}(F)$ defined by $(\mu M)(f):=\mu(M(f))$.

We denote by $\Vert M\Vert:=\sup_{f\in\mathcal{B}_{1}(E)}{\Vert M(f)\Vert} $
the norm of the operator $f\mapsto M(f)$ and we equip the Banach space
$\mathcal{M}(E)$ with the corresponding total variation norm $\Vert\mu
\Vert=\sup_{f\in\mathcal{B}_{1}(E)}|\mu(f)|$. We let $\beta(M)$ be the
Dobrushin coefficient of a bounded integral operator $M$ defined by the
following formula
\[
\beta(M):=\sup{\{\mbox{\rm osc}(M(f))\;;\;\;f\in\mbox{\rm Osc}_{1}(F)\}}%
\]
When $M$ has a constant mass, that is $M(1)(x)=M(1)(y)$ for any $(x,y)\in
E^{2}$, the operator $\mu\mapsto\mu M$ maps $\mathcal{M}_{0}(E)$ into
$\mathcal{M}_{0}(F)$ and $\beta(M)$ coincides with the norm of this operator.
We also denote by $(M^{k})_{k\geq0}$ the semigroup associated to $M$ given by
the recursive formulae $M^{k}(x,dz)=\int M^{k-1}(x,dy)M(y,dz)$, for $k\geq1$
and $M^{0}=Id$ the identity transition.

For any sequence of finite signed measures $\left(  \mu_{i}\right)  _{i\geq0}
$ defined on some collection of measurable spaces $\left(  E_{i}%
,\mathcal{E}_{i}\right)  _{i\geq0}$, we define $\left(  \mu_{i}\otimes\mu
_{j}\right)  \left(  dx_{i},dx_{j}\right)  =\mu_{i}\left(  dx_{i}\right)
\mu_{j}\left(  dx_{j}\right)  $, $\mu_{i}{}^{\otimes2}\left(  dx_{i}%
,dx_{i}^{\prime}\right)  =\left(  \mu_{i}\otimes\mu_{i}\right)  \left(
dx_{i},dx_{i}^{\prime}\right)  $ and
\[
\otimes_{i\leq k\leq j}\mu_{k}\text{ }\left(  dx_{i},dx_{i+1},\ldots
,dx_{j}\right)  =\prod_{k=i}^{j}\mu_{k}\left(  dx_{k}\right)  .
\]
For two bounded measurable functions $f$ and $g$ defined respectively on
$\left(  E,\mathcal{E}\right)  $ and $\left(  F,\mathcal{F}\right)  $ we
define $(f\otimes g)\left(  x,x^{\prime}\right)  =f\left(  x\right)  g\left(
x^{\prime}\right)  $.

We equip the set of distribution flows $\mathcal{M}(E)^{\mathbb{N}}$ with the
uniform total variation distance defined by
\[
\forall\eta=(\eta_{n})_{n\geq0},\text{ }\mu=(\mu_{n})_{n\geq0}\in
\mathcal{M}(E)^{\mathbb{N}}\qquad\Vert\eta-\mu\Vert:=\sup_{n\geq0}\Vert
\eta_{n}-\mu_{n}\Vert
\]
We extend a given bounded integral operator $\mu\in\mathcal{M}(E)\mapsto\mu
M\in\mathcal{M}(F)$ into an mapping
\[
\eta=(\eta_{n})_{n\geq0}\ \in\mathcal{M}(E)^{\mathbb{N}}\mapsto\eta
D=(\eta_{n}M)_{n\geq0}\in\mathcal{M}(F)^{\mathbb{N}}%
\]
Sometimes, we abuse the notation and denote by $\nu$ instead of $(\nu
)_{n\geq0}$ the constant distribution flow equal to a given measure $\nu
\in\mathcal{P}(E)$.

For any $\mathbb{R}^{d}$-valued function $f=(f^{i})_{1\leq i\leq d}%
\in\mathcal{B}(F)^{d}$, any integral operator $M$ from $E$ into $F$, and any
$\mu\in\mathcal{M}(F)$, we will slightly abuse the notation, and we write
$M(f)$ and $\mu(f)$ the $\mathbb{R}^{d}$-valued function and the point in
$\mathbb{R}^{d}$ given respectively by
\[
M(f):=\left(  M(f^{1}),\ldots,M(f^{d})\right)  \quad\mbox{\rm and}\quad
\mu(f):=\left(  \mu(f^{1}),\ldots,\mu(f^{d})\right)
\]
We also simplify the notation and sometimes we write
\[
M[(f^{1}-M(f^{1}))~(f^{2}-M(f^{2}))](x)
\]
instead of
\[%
\begin{array}
[c]{l}%
M[(f^{1}-M(f^{1})(x))~(f^{2}-M(f^{2})(x))](x)=M(f^{1}f^{2})(x)-M(f^{1}%
)(x)~M(f^{2})(x)
\end{array}
\]
Unless otherwise is stated, we denote by $c(k)$, $k\in\mathbb{N}$, a constant
whose value may vary from line to line but only depends on the parameter $k$.
Finally, we shall use $\sum_{\emptyset}=0$ and $\prod_{\emptyset}=1$.

\subsection{Examples\label{subsec:nonlinearexamples}}

We give here two classes of models where i-MCMC methods can be used.

\textit{Feynman-Kac models. }In this context, we have
\begin{equation}
\forall l\geq0\quad\forall(\mu,f)\in(\mathcal{P}(S^{(l)})\times\mathcal{B}%
(S^{(l+1)}))\qquad\Phi^{\left(  l+1\right)  }(\mu)(f):={\mu(G_{l}L_{l+1}%
(f))}/{\mu(G_{l})}\label{fksg}%
\end{equation}
where $L_{l+1}$ is a Markov transition kernel from $S^{(l)}$ into $S^{(l+1)}$
and $G_{l}:S^{(l)}\rightarrow\mathbb{R}^{+}$. In this situation, the solution
of the measure-valued equation (\ref{phi}) is given by%
\begin{equation}
\pi^{(l)}(f)={\gamma^{(l)}(f)}/{\gamma^{(l)}(1)}\quad\mbox{\rm with}\quad
\gamma^{(l)}(f):=\mathbb{E}\left(  f(X_{l})~\prod_{0\leq k<l}G_{k}%
(X_{k})\right) \label{eq:feynmankac}%
\end{equation}
where $(X_{l})_{l\geq0}$ is a Markov chain taking values in the state spaces
$(S^{(l)})_{l\geq0}$, with initial distribution $\pi^{(0)}$ and Markov
transitions $(L_{l})_{l\geq1}$. These Feynman-Kac models arise in a large
number of applications including nonlinear filtering, Bayesian statistics and
physics; see \cite{arnaud,fk}. Note that these models are quite flexible. For
instance, the reference Markov chain may represent the paths from the origin
up to the current time $l$ of an auxiliary Markov chain $\left(  X_{l}%
^{\prime}\right)  _{l\geq0}$ taking values in some state spaces $\left(
S_{l}^{\prime}\right)  _{l\geq0}$ with initial distribution $\pi^{(0)}$ and
Markov transitions $(L_{l}^{\prime})_{l\geq1}$; that is, we have
\begin{equation}
X_{l}:=(X_{0}^{\prime},\ldots,X_{l}^{\prime})\in S^{(l)}:=(S_{0}^{\prime
}\times\ldots\times S_{l}^{\prime})\label{pathmdels}%
\end{equation}
and consequently%
\begin{equation}
L_{l}\left(  x_{l-1},dy_{l}\right)  =\delta_{x_{l-1}}\left(  dy_{l-1}\right)
L_{l}^{\prime}\left(  y_{l-1}^{\prime},dy_{l}^{\prime}\right)
\label{eq:markovpotentialpathspaces}%
\end{equation}
When $G_{l}(x_{l})=G_{l}^{\prime}(x_{l}^{\prime})$, that is the potential
function only depends on the terminal value of the path $x_{l}=(x_{0}^{\prime
},\ldots,x_{l}^{\prime})$, the measures $\pi^{(l)}$ correspond to the path
space measures given for $l\geq1$ by
\begin{equation}
\pi^{(l)}\left(  dx_{l}\right)  \propto~\left\{  \prod_{0\leq k<l}%
G_{k}^{\prime}(x_{k}^{\prime})\right\}  ~\pi^{(0)}\left(  dx_{0}^{\prime
}\right)  \prod_{1\leq k\leq l}L_{k}^{\prime}\left(  x_{k-1}^{\prime}%
,dx_{k}^{\prime}\right) \label{eq:pathmeasures}%
\end{equation}
where `$\propto$' means `proportional to'. $\blacksquare$

\emph{Interacting annealing models}. These models were recently introduced in
\cite{atchade} and can be reinterpreted as a special case of (\ref{phi}). In
this scenario, we have $S^{(l)}=S$ and a pre-determined sequence of
probability distributions $(\pi^{(l)})_{l\geq0}$ of the form
\[
\pi^{\left(  l\right)  }(dx)=\frac{\exp{(-\beta}_{l}{V}\left(  x\right)
{)}~\lambda(dx)}{\lambda\left(  \exp{(-\beta}_{l}{V)}\right)  }%
\]
where $\lambda$ is a reference measure, $\left(  {\beta}_{l}\right)  _{l\geq
0}$ is an increasing positive sequence and $\lambda\left(  \exp{(-\beta}%
_{l}{V)}\right)  <\infty$. Based on this sequence, we build a sequence of
mappings $\left(  \Phi^{\left(  l\right)  }\right)  _{l\geq0}$ satisfying
(\ref{phi}) as follows. We introduce $\epsilon\in\lbrack0,1)$ and two
sequences of Markov kernels $\left(  K_{l}\right)  _{l\geq0}$ and $\left(
L_{l}\right)  _{l\geq0}$ where both $K_{l}$ and $L_{l}$ admit $\pi^{\left(
l\right)  }$ as invariant measure. We then set
\begin{equation}
\forall l\geq0\quad\forall(\mu,f)\in(\mathcal{P}(S^{(l)})\times\mathcal{B}%
(S^{(l+1)}))\qquad\Phi^{\left(  l+1\right)  }(\mu)\left(  f\right)  :=\Psi
_{l}(\mu)L_{l+1}K_{\epsilon,l+1}\left(  f\right) \label{eq:Atchademapping}%
\end{equation}
where $\Psi_{l}~:~\mu\in\mathcal{P}(S)\mapsto\Psi_{l}(\mu)\in\mathcal{P}(S)$
is defined by
\begin{equation}
\Psi_{l}(\mu)(dx):=\frac{G_{l}(x)~\mu(dx)}{\mu(G_{l})}\label{bgeq}%
\end{equation}
for $G_{l}\left(  x\right)  =\exp{(-(\beta_{l+1}-\beta_{l})V\left(  x\right)
)}$ and
\[
K_{\epsilon,l}:=(1-\epsilon)~\sum_{k\geq0}\epsilon^{k}~K_{l}^{k}%
\]
It is easy to check that (\ref{phi}) is satisfied for the mappings
(\ref{eq:Atchademapping}). $\blacksquare$

\subsection{Interacting Markov chain Monte Carlo methods}

\label{dmm}

We introduce a sequence of `initial' probability measures $\left(
\nu^{\left(  l\right)  }\right)  _{l\geq0}$ on $\left(  S^{(l)}\right)
_{l\geq0}$. We also introduce a Markov transition $M^{(0)}$ from $S^{(0)}$
into itself, and a collection of Markov transitions $M_{\mu}^{(l)}$ from
$S^{(l)} $ into itself, indexed by the parameter $l\geq0$, where $\mu
\in\mathcal{P}(S^{(l-1)})$. We further assume that the invariant measure of
each operator $M_{\mu}^{(l)}$ is given by $\Phi^{\left(  l\right)  }(\mu)$;
that is we have
\[
\forall l\geq0\quad\forall\mu\in\mathcal{P}(S^{(l-1)})\qquad\Phi^{\left(
l\right)  }(\mu)=\Phi^{\left(  l\right)  }(\mu)~M_{\mu}^{(l)}%
\]
For $l=0$, we use the convention $M_{\mu}^{(0)}=M^{(0)}$ and $\Phi^{\left(
0\right)  }(\mu)=\pi^{(0)}$.

We have now all the elements to define the i-MCMC algorithm. The algorithm
generates a sequence of processes $\left(  X^{(l)}\right)  _{l\geq0}$ where
$X^{(k)}:=(X_{n}^{(k)})_{n\geq0}$ is the process at level $k$ whose associated
occupation measure at iteration $n$ of the algorithm is denoted by%
\[
\eta_{n}^{(k)}:=\frac{1}{n+1}\sum_{p=0}^{n}\delta_{X_{p}^{(k)}}.
\]

At level $k=0$, $X^{(0)}$ is a Markov chain on $S^{(0)}$ with $X_{0}^{(0)}%
\sim\nu^{\left(  0\right)  }$ and Markov transitions $M^{(0)}$; that is
\[
\mathbb{P}(X_{n+1}^{(0)}\in dx~|~X_{n}^{(0)})=M^{(0)}\left(  X_{n}%
^{(0)},dx\right)  .
\]
At level $k\geq1$, given a realization of the chain $X^{(k-1)}$, the $k$-th
level chain $X^{(k)}$ is an inhomogeneous Markov chain with $X_{0}^{(k)}%
\sim\nu^{\left(  k\right)  }$ and Markov transitions $M_{\eta_{n}^{(k-1)}%
}^{(k)}$ at iteration $n$ depending on the current occupation measure
$\eta_{n}^{(k-1)}$ of the chain at level $(k-1)$; that is%
\begin{equation}
\mathbb{P}(X_{n+1}^{(k)}\in dx~|~X^{(k-1)},~X_{n}^{(k)})=M_{\eta_{n}^{(k-1)}%
}^{(k)}(X_{n}^{\left(  k\right)  },dx).\label{itdeff}%
\end{equation}
The rationale behind this is that the $k$-th level chain $X_{n}^{(k)}$ behaves
asymptotically as an homogeneous Markov chain with transition kernel
$M_{\pi^{(k-1)}}^{(k)}$ of invariant probability measure $\pi^{(k)}$ as long
as $\eta_{n}^{(k-1)}$ is a `good' approximation of $\pi^{(k-1)}$.

These i-MCMC\ algorithms can be interpreted as non-standard adaptive
MCMC\ schemes \cite{andrieumoulines2006,andrieuthoms2008,robertsrosenthal2007}
where the parameters to be adapted are probability measures instead of
finite-dimensional parameters. Algorithms relying on similar principles were
first proposed in \cite{andrieuajay} and independently in \cite{atchade}.
Related algorithms where we also have a sequence of nested MCMC-like chains
`feeding' each other have also recently appeared in statistics \cite{kou2006}
and physics \cite{lyman2006}.

We now give examples of such Markov kernels for the Feynman-Kac and
interacting annealing models described in section
\ref{subsec:nonlinearexamples}.

\textit{Feynman-Kac models. }Assume we are working on path spaces
$S^{(l)}:=(S^{(l-1)}\times S_{l}^{\prime})$ (cf. (\ref{pathmdels}%
)-(\ref{eq:pathmeasures})), we can select for $M_{\mu}^{(l)}$ a
Metropolis-Hastings kernel of independent proposal distribution $\left(
\mu\otimes L_{l}^{\prime}\right)  $ and target distribution $\Phi^{\left(
l\right)  }(\mu)$. More precisely, using the fact that
\[
\Phi^{\left(  l\right)  }(\mu)\left(  d(y_{l-1},y_{l}^{\prime})\right)
\propto\mu(dy_{l-1})~G_{l-1}^{\prime}(y_{l-1}^{\prime})~L_{l}^{\prime}%
(y_{l-1}^{\prime},dy_{l}^{\prime})
\]
the independent Metropolis-Hastings kernel $M_{\mu}^{(l)}$ using $\mu
(dy_{l-1})~L_{l}^{\prime}(y_{l-1}^{\prime},dy_{l}^{\prime})$ as a proposal
distribution is given by
\begin{equation}%
\begin{array}
[c]{l}%
M_{\mu}^{\left(  l\right)  }(x_{l},dy_{l})\\
\\
=\mu(dy_{l-1})~L_{l}^{\prime}(y_{l-1}^{\prime},dy_{l}^{\prime})~\left(
1\wedge\frac{G_{l-1}^{\prime}(y_{l-1}^{\prime})}{G_{l-1}^{\prime}%
(x_{l-1}^{\prime})}\right)  +\left(  1-\mu\left(  1\wedge\frac{G_{l-1}%
^{\prime}(y_{l-1}^{\prime})}{G_{l-1}^{\prime}(x_{l-1}^{\prime})}\right)
\right)  ~\delta_{x_{l}}(dy_{l})
\end{array}
\label{eq:metropolishastingsfeynmankac}%
\end{equation}
where we recall that $x_{l}=\left(  x_{l-2},x_{l-1}^{\prime},x_{l}^{\prime
}\right)  =\left(  x_{l-1},x_{l}^{\prime}\right)  \in S^{(l)}=(S^{(l-1)}\times
S_{l}^{\prime})$ and $y_{l}=\left(  y_{l-1},y_{l}^{\prime}\right)  \in
S^{(l)}=(S^{(l-1)}\times S_{l}^{\prime})$. $\blacksquare$

\emph{Interacting annealing models}. In this case, we can select
\begin{equation}
M_{\mu}^{(l)}(x,dy)=\epsilon K_{l}(x,dy)+(1-\epsilon)~\Psi_{l-1}(\mu
)L_{l}(dy).\label{ouf}%
\end{equation}
One can easily check that $M_{\mu}^{(l)}$ admits $\Phi^{\left(  l\right)
}(\mu)$ as invariant probability measure. $\blacksquare$

For sufficiently regular models, we proved in~\cite{brockwell,dd08} that the
occupations measures $\eta_{n}^{(l)}$ converge to the solution $\pi^{(l)}$ of
the equation (\ref{phi}), in the sense that $\lim_{n\rightarrow\infty}\eta
_{n}^{(l)}(f)=\pi^{(l)}(f)$ almost surely for $f\in\mathcal{B}(S^{(l)})$. The
articles~\cite{bdd08,dd08} also provide a collection of non asymptotic
$\mathbb{L}_{r}$-mean error estimates and exponential deviations inequalities.
The fluctuation analysis of $\eta_{n}^{(l)}$ around the limiting measure
$\pi^{(l)}$ has been initiated in~\cite{bdd08} in the special case where
$M_{\mu}^{(l)}(x_{l},\mbox{\LARGE .})=\Phi^{\left(  l\right)  }(\mu)$. In this
`simpler' situation, the $l$-th level chain $X^{(l)}=\left(  X_{n}%
^{(l)}\right)  $ is given $X^{(l-1)}$ a collection of conditionally
independent random variables with $X_{0}^{(l)}\sim\nu^{\left(  l\right)  }$
and $X_{n}^{(l)}\sim\Phi^{\left(  l\right)  }(\eta_{n-1}^{(l-1)})$ for
$n\geq1$.

\subsection{Contribution and organization of the paper}

The present article studies the fluctuations of the occupation measures
$\left(  \eta_{n}^{(l)}\right)  _{l\geq0}$ associated to the class of i-MCMC
algorithms towards their limiting values $\left(  \pi^{(l)}\right)  _{l\geq0}
$. Briefly speaking, our analysis proceeds as follows. First, we study
weighted sequences of local random fields $V_{n}^{(l)}$ which are related to
the fluctuations of the occupation measures $\eta_{p}^{(l)}$ around their
local invariant measures $\Phi^{\left(  l\right)  }\left(  \eta_{p-1}%
^{(l-1)}\right)  $ for $p\leq n$. We show that these random fields
$(V_{n}^{(l)})_{l\geq0}$ converges in law, as $n$ tends to infinity and in the
sense of finite dimensional distributions, to a sequence of independent and
centered Gaussian fields $\left(  V^{(l)}\right)  _{l\geq0}$ with covariance
functions defined in terms of the resolvent operator associated to the Markov
transition $M_{\pi^{(l-1)}}^{(l)}$ and its invariant probability measure
$\pi^{(l)}$. Finally, we deduce the fluctuations of $\eta_{n}^{(l)}$ around
their limiting values $\pi^{(l)}$ by a simple application of the continuous
mapping theorem (or the multivariate $\delta$-method) applied to a first order
decomposition of the error $\sqrt{n}\left[  \eta_{n}^{(l)}-\pi^{(l)}\right]  $
in terms of the random fields $(V_{n}^{(k)})_{0\leq k\leq l}$.

The rest of the paper is organized as follows. The main result of the article
is presented in full details in section~\ref{ssr}. The regularity conditions
are summarized in section~\ref{reg}. In section~\ref{ftcl} we state a
multivariate CLT in terms of the semigroup associated with a first order
expansion of the mappings $\Phi^{\left(  l\right)  }$ appearing in
(\ref{phi}). Section~\ref{fnh} addresses the fluctuation analysis of an
abstract class of time inhomogeneous Markov chains. In section~\ref{resl}, we
present a preliminary resolvent analysis to estimate the regularity properties
of resolvent operators and invariant measure type mappings. In
section~\ref{lln}, we apply these results to study the local fluctuations of a
class of weighted occupation measures associated to self--interacting chains.
Section~\ref{secfrf} addresses the fluctuation analysis of local interaction
random fields associated with i-MCMC algorithms. The proof of the main theorem
presented in section~\ref{ftcl} is a direct consequence of a fluctuation
theorem for local interaction random fields, and it is given at the end of
section~\ref{introf}. Finally, we establish in section
\ref{sec:pathspacemodels} that the regularity conditions discussed in
section~\ref{reg} are also valid for a path space extension of i-MCMC\ algorithms.

\section{Statement of some results}

\label{ssr}

\subsection{Regularity conditions\label{sec:regularityconditions}}

\label{reg}Our first regularity condition is a first order weak regularity
condition on the mappings $\Phi^{\left(  l\right)  }$ governing the
measure-valued equation (\ref{phi}). We assume that, for any $l\geq0$, the
mappings $\Phi^{\left(  l+1\right)  }:\mathcal{P}(S^{(l)})\rightarrow
\mathcal{P}(S^{(l+1)})$ satisfy the following first order local decomposition
\begin{equation}
\left[  \Phi^{\left(  l+1\right)  }(\mu)-\Phi^{\left(  l+1\right)  }%
(\eta)\right]  =(\mu-\eta)D_{l+1}+\Xi_{l}(\mu,\eta)\label{firsto}%
\end{equation}
where $D_{l+1}~:~\mathcal{B}(S^{(l+1)})\rightarrow\mathcal{B}(S^{(l)})$ is a
bounded integral operator that may depend on the measure $\eta$ and $\Xi
_{l}(\mu,\eta)$ is a remainder signed measure on $S^{(l+1)}$ indexed by the
set of probability measures $\mu,\eta\in\mathcal{P}(S^{(l)})$. We further
require that
\begin{equation}
\left\vert \Xi_{l}(\mu,\eta)(f)\right\vert \leq\int~\left\vert (\mu
-\eta)^{\otimes2}(g)\right\vert ~\varXi_{l}(f,dg)\label{condxi}%
\end{equation}
for some integral operator $\varXi_{l}$ from $\mathcal{B}(S^{(l+1)})$ into the
set $\mathcal{T}_{2}(S^{(l)})$ of all tensor product functions $g=\sum_{i\in
I}~a_{i}~(h_{i}^{1}\otimes h_{i}^{2})$, with $I\subset\mathbb{N}$, $(h_{i}%
^{1},h_{i}^{2})_{i\in I}\in(\mathcal{B}(S^{(l)})^{2})^{I}$, and a sequence of
numbers $(a_{i})_{i\in I}\in\mathbb{R}^{I}$ such that
\begin{equation}
|g|:=\sum_{i\in I}~|a_{i}|~\Vert h_{i}^{1}\Vert\Vert h_{i}^{2}\Vert
<\infty\quad\mbox{\rm and}\quad\chi_{l}:=\sup_{f\in\mathcal{B}_{1}(S^{(l+1)}%
)}\int~|g|~\varXi_{l}(f,dg)<\infty\label{defchi}%
\end{equation}

Our second set of regularity conditions are for the Markov kernels $M_{\mu
}^{(l)}$. We assume these kernels satisfy the following two regularity
conditions
\begin{equation}
m_{l}(n_{l}):=\sup_{\mu\in\mathcal{P}(S^{(l-1)})}{\beta((M_{\mu}^{(l)}%
)^{n_{l}})}<1\label{lipcna}%
\end{equation}
and
\begin{equation}
\left\Vert \left[  M_{\mu}^{(l)}-M_{\nu}^{(l)}\right]  (f)\right\Vert \leq
\int~\left\vert \left[  \mu-\nu\right]  (g)\right\vert ~\Gamma_{l,\mu
}(f,dg)\label{lipcnb}%
\end{equation}
for some collection of bounded integral operators $\Gamma_{l,\mu}$ from
$\mathcal{B}(S^{(l)})$ into $\mathcal{B}(S^{(l-1)})$ and indexed by the set of
measures $\mu\in\mathcal{P}(S^{(l-1)})$ with
\[
\sup_{\mu\in\mathcal{P}(S^{(l-1)})}\int~\Gamma_{l,\mu}(f,dg)~\Vert g\Vert
\leq\Lambda_{l}~\Vert f\Vert\quad\mbox{\rm and}\quad\Lambda_{l}<\infty
\]

We end this section with some comments on this set of conditions.

The regularity condition (\ref{firsto})-(\ref{condxi}) is a first order
refinement of a Lipschitz type condition we used in ~\cite{brockwell,dd08} to
derive a series of $\mathbb{L}_{p}$-mean error bounds and exponential
inequalities. This condition has been introduced in~\cite{bdd08} for studying
the fluctuations of the simple i-MCMC algorithm corresponding to $M_{\mu
}^{(l)}(x,\mbox{\LARGE .})=\Phi^{\left(  l\right)  }(\mu)$. The regularity
condition (\ref{lipcna}) is an ergodicity condition on the Markov kernel
${M_{\mu}^{(l)}}$. Finally the regularity condition (\ref{lipcnb}) is a local
Lipschitz type continuity condition on the kernel $M_{\mu}^{(l)}$. This
condition is less stringent than the one used in~\cite{dd08} where it is
assumed that (\ref{lipcnb}) holds for some operators $\Gamma_{l,\mu}%
=\Gamma_{l}$ that do not depend on $\mu$. Therefore, most of the asymptotic
results presented in~\cite{dd08} do not apply in the present context.
Nevertheless, it can be checked that the inductive proof of the $\mathbb{L}%
_{p}$-mean error bounds presented in theorem 5.2 in \cite{bdd08} hold true
under the weaker condition (\ref{lipcnb}); thus, for every $l\geq0 $ and any
function $f\in\mathcal{B}(S^{(l)})$, we know that $\eta_{n}^{(l)}(f)$
converges almost surely to $\pi^{(l)}(f)$ as $n\rightarrow\infty$. The main
advantage of the set of conditions presented here is that it is stable under a
state space enlargement (see Section \ref{sec:pathspacemodels}), so that the
asymptotic analysis of such algorithms, including the multivariate CLT
presented in the next section, applies directly without further work to i-MCMC
algorithms on path spaces.

We illustrate these regularity conditions for the models discussed in
section~\ref{subsec:nonlinearexamples}. We further assume in the rest of this
section that $\left(  G_{l}\right)  _{l\geq0}$ is a collection of $\left(
0,1\right]  $-valued potential functions on some state space $\left(
S^{(l)}\right)  _{l\geq0}$ such that
\begin{equation}
\forall l\geq0\qquad\inf_{S^{(l)}}{G_{l}}>0\label{eq:conditionpotential}%
\end{equation}

\textit{Feynman-Kac models. }To establish (\ref{firsto})-(\ref{condxi}), we
observe that the mapping $\Psi_{l}$ defined in (\ref{bgeq}) can be rewritten
in terms of a nonlinear transport equation
\[
\Psi_{l}(\mu)(dy)=\left(  \mu\mathcal{S}_{l,\mu}\right)  (dy):=\int%
\mu(dx)\mathcal{S}_{l,\mu}(x,dy)
\]
where
\[
\mathcal{S}_{l,\mu}(x,dy)=G_{l}(x)~\delta_{x}(dy)+\left(  1-G_{l}(x)\right)
~\Psi_{l}(\mu)(dy)
\]
Using the decomposition
\begin{equation}
\Psi_{l}(\mu)-\Psi_{l}(\eta)=(\mu-\eta)\mathcal{S}_{l,\eta}+\mu(\mathcal{S}%
_{l,\mu}-\mathcal{S}_{l,\eta})\Rightarrow\Psi_{l}(\mu)-\Psi_{l}(\eta)=\frac
{1}{\mu(G_{l})}~(\mu-\eta)\mathcal{S}_{l,\eta}\label{decouf}%
\end{equation}
we prove the first order decomposition
\[
\Psi_{l}(\mu)-\Psi_{l}(\eta)=(\mu-\eta)D_{l}^{\prime}+\Xi_{l-1}^{\prime}%
(\mu,\eta)
\]
with the integral operators $D_{l}^{\prime}$ defined for any $f\in
\mathcal{B}(S^{(l)})$ by $D_{l}^{\prime}(f):=(\eta(G_{l}))^{-1}~\mathcal{S}%
_{l,\eta}(f)$, and the remainder measures
\[
\Xi_{l-1}^{\prime}(\mu,\eta)(f):=\left[  \frac{1}{\mu(G_{l})}-\frac{1}%
{\eta(G_{l})}\right]  ~(\mu-\eta)\mathcal{S}_{l,\eta}(f)
\]
Using the fact that
\[
\left[  \frac{1}{\mu(G_{l})}-\frac{1}{\eta(G_{l})}\right]  =\frac{(\eta
-\mu)(G_{l})}{\mu(G_{l})\eta(G_{l})}%
\]
and
\[
(\mu-\eta)\mathcal{S}_{l,\eta}(f)=(\mu-\eta)\left(  G_{l}f\right)  -(\mu
-\eta)(G_{l})~\Psi_{l}(\eta)(f)
\]
we obtain%
\begin{align*}
\left\vert \Xi_{l-1}^{\prime}(\mu,\eta)(f)\right\vert  & \leq\frac{1}{\inf
G_{l}^{2}}\left[  \left\vert (\mu-\eta)^{\otimes2}(G_{l}\otimes(G_{l}%
f))\right\vert +\Vert f\Vert\left\vert (\mu-\eta)^{\otimes2}(G_{l}\otimes
G_{l})\right\vert \right] \\
:=  & \int~\left\vert (\mu-\eta)^{\otimes2}(g)\right\vert ~\varXi_{l-1}%
^{\prime}(f,dg)
\end{align*}
with the integral operator
\[
\varXi_{l-1}^{\prime}(f,dg)=\frac{1}{\inf G_{l}^{2}}~\left(  \delta
_{G_{l}\otimes(G_{l}f)}(dg)+\Vert f\Vert~\delta_{G_{l}\otimes G_{l}%
}(dg)\right)
\]
We check that the mappings (\ref{fksg}) satisfy (\ref{firsto}) using the fact
that
\begin{equation}
\Phi^{\left(  l+1\right)  }(\mu)-\Phi^{\left(  l+1\right)  }(\eta)=(\Psi
_{l}(\mu)-\Psi_{l}(\eta))L_{l+1}=(\mu-\eta)D_{l+1}+\Xi_{l}(\mu,\eta
)\label{eq:cond1FK}%
\end{equation}
with the first order operator $D_{l+1}=D_{l}^{\prime}L_{l+1}$ and the
remainder measure $\Xi_{l}(\mu,\eta)=\Xi_{l-1}^{\prime}(\mu,\eta)L_{l+1}$. The
remainder measure satisfies (\ref{condxi}) for
\begin{equation}
\varXi_{l}(f,dg)=\frac{1}{\inf G_{l}^{2}}~\left(  \delta_{G_{l}\otimes
(G_{l}L_{l+1}(f))}(dg)+\Vert L_{l+1}(f)\Vert~\delta_{G_{l}\otimes G_{l}%
}(dg)\right) \label{eq:cond2FK}%
\end{equation}
We mention that in this case the parameter $\chi_{l}$ defined in
(\ref{defchi}) is such that
\[
\chi_{l}\leq2\sup G_{l}^{2}/\inf G_{l}^{2}%
\]

Assume we are working on path spaces $S^{(l)}:=(S^{(l-1)}\times S_{l}^{\prime
})$ where $G_{l}\left(  x_{l}\right)  =G_{l}^{\prime}\left(  x_{l}^{\prime
}\right)  $ and $\left(  S_{l}^{\prime}\right)  _{l\geq0}$ are finite spaces.
If we use for $M_{\mu}^{(l)}\left(  x,dy\right)  $ the independent
Metropolis-Hastings kernel (\ref{eq:metropolishastingsfeynmankac}),
(\ref{lipcna}) is satisfied as $\left\Vert G_{l}^{\prime}\right\Vert \leq1$;
e.g. \cite[theorem 2.1.]{mengersentweedie1996}. Additionally, we have
\[%
\begin{array}
[c]{l}%
\left\vert \int\left(  M_{\mu}^{\left(  l\right)  }(x_{l},dy_{l})-M_{\nu
}^{(l)}(x_{l},dy_{l})\right)  f\left(  y_{l}\right)  \right\vert \\
\leq\left\vert \int\left(  \mu-\nu\right)  (dy_{l-1})~L_{l}^{\prime}%
(y_{l-1}^{\prime},dy_{l}^{\prime})~\left(  1\wedge\frac{G_{l-1}^{\prime
}(y_{l-1}^{\prime})}{G_{l-1}^{\prime}(x_{l-1}^{\prime})}\right)  f\left(
y_{l}\right)  \right\vert +\left\Vert f\right\Vert \left\vert \int\left(
\mu-\nu\right)  (dy_{l-1})\left(  1\wedge\frac{G_{l-1}^{\prime}(y_{l-1}%
^{\prime})}{G_{l-1}^{\prime}(x_{l-1}^{\prime})}\right)  \right\vert .
\end{array}
\]
So (\ref{lipcnb}) is satisfied for
\[
\Gamma_{l,\mu}(f,dg)=\sum_{x_{l-1}^{\prime}\in S_{l-1}^{\prime}}\delta
_{L_{l}^{\prime}(f)~\left(  1\wedge\frac{G_{l-1}^{\prime}}{G_{l-1}^{\prime
}(x_{l-1}^{\prime})}\right)  }\left(  g\right)  +\left\Vert f\right\Vert
\delta_{1\wedge\frac{G_{l-1}^{\prime}}{G_{l-1}^{\prime}(x_{l-1}^{\prime})}%
}\left(  g\right)
\]
where%
\begin{align*}
\int~\Gamma_{l,\mu}(f,dg)~\Vert g\Vert & =\sum_{x_{l-1}^{\prime}\in
S_{l-1}^{\prime}}~\Vert L_{l}^{\prime}(f)~\left(  1\wedge\frac{G_{l-1}%
^{\prime}}{G_{l-1}^{\prime}(x_{l-1}^{\prime})}\right)  \Vert+\left\Vert
f\right\Vert \left\Vert 1\wedge\frac{G_{l-1}^{\prime}}{G_{l-1}^{\prime
}(x_{l-1}^{\prime})}\right\Vert \\
& \leq2\left(  \sum_{x_{l-1}^{\prime}\in S_{l-1}^{\prime}}\left\Vert
1\wedge\frac{G_{l-1}^{\prime}}{G_{l-1}^{\prime}(x_{l-1}^{\prime})}\right\Vert
\right)  \text{ }\left\Vert f\right\Vert .
\end{align*}
$\blacksquare$

\emph{Interacting annealing models}. To establish (\ref{firsto})-(\ref{condxi}%
), we can proceed similarly to Feynman-Kac models. It is sufficient to
substitute $L_{l+1}K_{\epsilon,l+1}$ to $L_{l+1}$ in (\ref{eq:cond1FK}%
)-(\ref{eq:cond2FK}). In this context, the Markov transitions given in
(\ref{ouf}) are such that
\[
M_{\eta}^{(l)}(x,dy)\geq(1-\epsilon)~\Psi_{l-1}(\mu)L_{l}(dy)\Longrightarrow
\beta(M_{\eta}^{(l)})\leq\epsilon.
\]
so (\ref{lipcna}) is satisfied with $n_{l}=1$ and $m_{l}(1)\leq\epsilon$.
Moreover, we have
\[
\left[  M_{\mu}^{(l)}-M_{\nu}^{(l)}\right]  (f)=(1-\epsilon)~\left[
\Psi_{l-1}(\mu)-\Psi_{l-1}(\nu)\right]  L_{l}(f)
\]
Using the decomposition (\ref{decouf}) one proves that
\[
\left\Vert \left[  M_{\mu}^{(l)}-M_{\nu}^{(l)}\right]  (f)\right\Vert
\leq\frac{1-\epsilon}{\inf G_{l-1}}~\left\vert (\mu-\nu)\left(  \mathcal{S}%
_{l-1,\mu}L_{l}(f)\right)  \right\vert
\]
Hence condition (\ref{lipcnb}) is satisfied with $\Gamma_{l,\mu}%
(f,dg)=\frac{1-\epsilon}{\inf G_{l-1}}~\delta_{\mathcal{S}_{l-1,\mu}L_{l}%
(f)}(dg) $ and $\Lambda_{l}\leq{(1-\epsilon)}/{\inf G_{l-1}}$. $\blacksquare$

\subsection{A multivariate central limit theorem}

\label{ftcl}

To describe precisely the fluctuations of the empirical measures $\eta
_{n}^{(l)}$ around their limiting value $\pi^{(l)}$, we need a few additional
notations. We denote by $D_{k,l}$ with $0\leq k\leq l$ the semigroup
associated with the bounded integral operators $D_{k}$ introduced in
(\ref{firsto}). More formally, we have
\[
\forall1\leq k\leq l\qquad D_{k,l}=D_{k}D_{k+1}\ldots D_{l}%
\]
For $k>l$, we use the convention $D_{k,l}=Id$, the identity operator.

Using this notation, the multivariate CLT describing the fluctuations of the
i-MCMC algorithm around the solution of the equation (\ref{phi}) is stated as
follows. We remind the reader that the integral operator $D_{l+1}$ may depend
on the measure $\pi^{\left(  l\right)  }$.

\begin{theo}
\label{tcl} For every $k\geq0$, the sequence of random fields $(U_{n}%
^{(k)})_{n\geq0}$ on $\mathcal{B}(S^{(k)})$ defined below
\[
U_{n}^{(k)}:=\sqrt{n+1}~\left[  \eta_{n}^{(k)}-\pi^{(k)}\right]
\]
converges in law, as $n$ tends to infinity and in the sense of finite
dimensional distributions, to a sequence of Gaussian random fields $U^{(k)}$
on $\mathcal{B}(S^{(k)})$ given by the following formula
\begin{equation}
U^{(k)}:=\sum_{0\leq l\leq k}~\frac{\sqrt{(2l)!}}{l!}~~V^{(k-l)}%
D_{(k-l)+1,k}\label{eq:asymptoticvariance}%
\end{equation}
Here $\left(  V^{(l)}\right)  _{l\geq0}$ stands for a collection of
independent and centered Gaussian fields with a variance function given by
\begin{equation}%
\begin{array}
[c]{l}%
\mathbb{E}\left(  V^{(l)}(f)^{2}\right) \\
\\
=\pi^{(l)}\left[  (f-\pi^{(l)}(f))^{2}\right]  +2\sum_{n\geq1}\pi^{(l)}\left[
(f-\pi^{(l)}(f))~\left(  M_{\pi^{(l-1)}}^{(l)}\right)  ^{n}(f-\pi
^{(l)}(f))\right]
\end{array}
\label{cov}%
\end{equation}

\end{theo}

In the special case where $M_{\mu}^{(l)}(x,\mbox{\LARGE .})=\Phi^{\left(
l\right)  }(\mu)$ for all $l\geq1$, that is $\mathbb{E}\left(  V^{(l)}%
(f)^{2}\right)  =\pi^{(l)}\left[  (f-\pi^{(l)}(f))^{2}\right]  $, the result
corresponds to the one obtained\ previously in~\cite{bdd08}. This special
class of i-MCMC algorithms behaves as a sequence of independent random
variables with distributions $\Phi^{\left(  l\right)  }(\eta_{n}^{(l-1)})$
given by the local invariant measures of MCMC\ chains with transition kernels
$M_{\eta_{n}^{(l-1)}}^{(l)}(x,\mbox{\LARGE .})$. In the more general case
considered here, the additional terms on the right hand side of (\ref{cov})
reflects the fluctuations of these MCMC algorithms around their limiting
invariant probability measures.

Finally we note that, if it was possible to sample exactly from $\pi^{\left(
k-1\right)  }$, then we would have $U^{(k)}=V^{\left(  k\right)  }$. However,
we need to approximate $\pi^{\left(  k-1\right)  }$ using an MCMC\ kernel
which itself relies on an MCMC\ approximation of $\pi^{\left(  k-2\right)  }$
and so on. The price to pay for these additional approximations appears
clearly in (\ref{eq:asymptoticvariance}).

\textit{A Toy Example.} Consider a Feynman-Kac model where\ $S^{(l)}%
:=(S^{(l-1)}\times S)$ with $S=\left\{  1,2\right\}  $ and
\[
G_{l}^{\prime}\left(  1\right)  =p^{\left(  \beta_{l+1}-\beta_{l}\right)
},\text{ }G_{l}^{\prime}\left(  2\right)  =q^{\left(  \beta_{l+1}-\beta
_{l}\right)  }%
\]
for $p=1-q>0$, $\left(  {\beta}_{l}\right)  _{l\geq0}$ is an increasing
positive sequence and
\[
L_{l+1}^{\prime}=\left(
\begin{tabular}
[c]{ll}%
$1-\pi_{l+1}\left(  2\right)  $ & $\pi_{l+1}\left(  2\right)  $\\
$\pi_{l+1}\left(  1\right)  $ & $1-\pi_{l+1}\left(  1\right)  $%
\end{tabular}
\right)
\]
where
\[
\pi_{l+1}\left(  1\right)  =1-\pi_{l+1}\left(  2\right)  =\frac{p^{\beta
_{l+1}}}{p^{\beta_{l+1}}+q^{\beta_{l+1}}}.
\]
We have $\pi_{l+1}L_{l+1}^{\prime}=\pi_{l+1}$ with $\pi_{l+1}=\left(
\pi_{l+1}\left(  1\right)  \text{ }\pi_{l+1}\left(  2\right)  \right)  $ and
it is easy to check that the resulting Feynmac-Kac probability measure
$\pi^{\left(  l+1\right)  }$ on $S^{l+2}$ admits as a marginal distribution
$\pi_{l+1}$ at time $l+1.$ For sake of illustration, we derive the expression
of all the terms appearing in the variance of $U^{(k)}\left(  f\right)  $
in\ equation (\ref{eq:asymptoticvariance}) of Theorem \ref{tcl}.

In this scenario, we have established in section
\ref{sec:regularityconditions} that $D_{l+1}=D_{l}^{\prime}L_{l+1}$ with
$D_{l}^{\prime}=S_{l,\pi^{\left(  l\right)  }}/\pi^{\left(  l\right)  }\left(
G_{l}\right)  =S_{l,\pi^{\left(  l\right)  }}/\pi_{l}\left(  G_{l}^{\prime
}\right)  $ where, recalling the notation $x_{l}=\left(  x_{0}^{\prime}%
,\ldots,x_{l}^{\prime}\right)  $, we obtain%
\begin{align*}
S_{l,\pi^{\left(  l\right)  }}\left(  x_{l},y_{l}\right)   & =G_{l}^{\prime
}\left(  x_{l}^{\prime}\right)  \delta_{x_{l}}\left(  y_{l}\right)  +\left(
1-G_{l}^{\prime}\left(  x_{l}^{\prime}\right)  \right)  \frac{\pi^{\left(
l\right)  }\left(  y_{l}\right)  G_{l}^{\prime}\left(  y_{l}^{\prime}\right)
}{\pi_{l}\left(  G_{l}^{\prime}\right)  },\\
L_{l+1}\left(  x_{l},y_{l+1}\right)   & =\delta_{x_{l}}\left(  y_{l}\right)
L_{l+1}^{\prime}\left(  y_{l}^{\prime},y_{l+1}^{\prime}\right)
\end{align*}
so%
\[
D_{l+1}\left(  x_{l},y_{l+1}\right)  =\frac{1}{\pi_{l}\left(  G_{l}^{\prime
}\right)  }\left\{  G_{l}^{\prime}\left(  x_{l}^{\prime}\right)  \delta
_{x_{l}}\left(  y_{l}\right)  L_{l+1}^{\prime}\left(  y_{l}^{\prime}%
,y_{l+1}^{\prime}\right)  +\left(  1-G_{l}^{\prime}\left(  x_{l}^{\prime
}\right)  \right)  \pi^{\left(  l+1\right)  }\left(  y_{l+1}\right)  \right\}
.
\]
Finally the Markov transition kernel $M_{\pi^{\left(  l\right)  }}^{\left(
l\right)  }(x_{l},y_{l})$ satisfies
\[
M_{\pi^{\left(  l\right)  }}^{\left(  l\right)  }(x_{l},y_{l})=\pi^{\left(
l\right)  }(y_{l-1})~L_{l}^{\prime}(y_{l-1}^{\prime},y_{l}^{\prime})~\left(
1\wedge\frac{G_{l-1}^{\prime}(y_{l-1}^{\prime})}{G_{l-1}^{\prime}%
(x_{l-1}^{\prime})}\right)  +\left(  1-\pi^{\left(  l\right)  }\left(
1\wedge\frac{G_{l-1}^{\prime}(y_{l-1}^{\prime})}{G_{l-1}^{\prime}%
(x_{l-1}^{\prime})}\right)  \right)  ~\delta_{x_{l}}(y_{l}).
\]
Clearly even in this toy example, the expression of the variance of
$U^{(k)}\left(  f\right)  $ is unfortunately analytically intractable. It is
additionally only possible to compute numerically this variance for very small
values of $k$ as the semigroup $D_{1,k}$ and the transition matrix
$M_{\pi^{\left(  k\right)  }}^{\left(  k\right)  }$ would have to be computed
for a number of values increasing exponentially fast with $k$. $\blacksquare$

Since the original version of this paper \cite{bdd08b}, the authors have
become aware of a recent paper of Yves Atchadé~\cite{atchade} which analyzes
the local fluctuations of some related algorithms; namely the
importance-resampling MCMC algorithm which is an interacting annealing model
as described in section \ref{subsec:nonlinearexamples} and a version of the
equi-energy sampler \cite{kou2006}. In \cite{atchade}, the author provides a
CLT associated to some random measures $\pi_{2}^{(k)}$ which converge almost
surely to $\pi_{2}^{(k)}$ as $n$ tends to infinity. We provide here a
multivariate CLT\ valid for any $l$ (and not only $l=2$). The variance
expression we obtain depends explicitly on the first order semigroup $D_{k,l}
$ and the fluctuations of the i-MCMC algorithm on lower indexed levels.

\section{On the fluctuations of time inhomogeneous Markov chains}

\label{fnh}

To establish the proof of our main result (theorem \ref{tcl}), it is necessary
to first provide some general results about the fluctuations of time
inhomogeneous Markov chains with Markov transitions that may depend on some
predictable flow of distributions on some possibly different state space.

\subsection{Description of the model}

\label{dmmc} We consider a collection of Markov transitions $M_{\eta}$ on some
measurable space $(S,\mathcal{S})$ indexed by the set of probability measures
$\eta\in\mathcal{P}(S^{\prime})$, on some possibly different measurable space
$(S^{\prime},\mathcal{S}^{\prime})$. We further assume that there exists an
integer $n_{0}\geq0$ such that
\begin{equation}
m(n_{0}):=\sup_{\eta\in\mathcal{P}(S^{\prime})}{\beta(M_{\eta}^{n_{0}}%
)}<1\quad\mbox{\rm and we set}\quad p(n_{0}):={2n_{0}}/{(1-m(n_{0}%
))}\label{lipcnp1}%
\end{equation}
We also assume that for any pair of measures $(\eta,\mu)\in\mathcal{P}%
(S^{\prime})^{2}$ we have
\begin{equation}
\left\Vert \left[  M_{\mu}-M_{\eta}\right]  (f)\right\Vert \leq\int~\left\vert
\left[  \mu-\eta\right]  (g)\right\vert ~\Gamma_{\mu}(f,dg)\label{lipcnp2}%
\end{equation}
for some collection of bounded integral operator $\Gamma_{\mu}$ from
$\mathcal{B}(S)$ into $\mathcal{B}(S^{\prime})$, indexed by the set of
measures $\mu\in\mathcal{P}(S^{\prime})$ with
\[
\sup_{\mu\in\mathcal{P}(S^{\prime})}\int~\Gamma_{\mu}(f,dg)~\Vert g\Vert
\leq\Lambda~\Vert f\Vert\quad\mbox{\rm for some finite constant}\quad
\Lambda<\infty
\]
We consider an increasing sequence of $\sigma$-fields $\left(  \mathcal{F}%
_{n}\right)  _{n\geq0}$ on some probability space $\left(  \Omega
,\mathcal{F},\mathbb{P}\right)  $. We let $\eta_{n}$ be a $\mathcal{P}%
(S^{\prime})$-valued random process adapted to the filtration $\mathcal{F}%
_{n}$ (i.e. each probability distribution $\eta_{n}$ is $\mathcal{F}_{n}%
$-measurable). We further assume that $\mathcal{F}_{n}$ contains the $\sigma
$-field generated by the random states $X_{p}$ from the origin $p=0$ up to the
current time horizon $p=n$ of an $S$-valued non homogeneous Markov chain
$X_{n}$ with a prescribed initial distribution $\nu\in\mathcal{P}(S)$, and
some transitions defined by
\begin{equation}
\forall n\geq0\qquad\mathbb{P}(X_{n+1}\in dx~|~\mathcal{F}_{n})=M_{\eta_{n}%
}(X_{n},dx)\label{eqqM}%
\end{equation}
For example, $\mathcal{F}_{n}$ could be $\sigma\left(  (X_{p}^{\prime}%
,X_{p}),~0\leq p\leq n\right)  $, i.e. the canonical sigma field associated
with the $(S^{\prime}\times S)$-valued process $(X_{n}^{\prime},X_{n}%
)_{n\geq0}$, and $\eta_{n}=\frac{1}{(n+1)}\sum_{p=0}^{n}\delta_{X_{p}^{\prime
}}$ is the flow of occupation measures of $\left(  X_{n}^{\prime}\right)
_{n\geq0}$. In this context, (\ref{eqqM}) reflects the fact that, given
$(\eta_{n})_{n\geq0}$, the process $\left(  X_{n}\right)  _{n\geq0}$ is a
Markov chain with random Markov transitions defined in terms of the occupation
measures $\left(  \eta_{n}\right)  _{n\geq0}$.

We further assume that the variations of $\left(  \eta_{n}\right)  _{n}$ are
controlled by some sequence of random variables $\tau(n)$ in the sense that
\begin{equation}
\forall n\geq0\qquad\Vert\eta_{n}-\eta_{n-1}\Vert\leq\tau(n),\quad\mbox{\rm
and we set}\quad\overline{\tau}(n):=\sum_{0\leq p\leq n}\tau(p)\label{lipcnp3}%
\end{equation}
For $n=0$ we use the convention $\eta_{-1}=0$, the null measure on $S^{\prime
}$.

\subsection{Regularity properties of resolvent operators}

\label{resl}

The main simplification of conditions (\ref{lipcnp1}) comes from the fact that
$M_{\eta}$ has an unique invariant measure
\[
\Phi(\eta)M_{\eta}=\Phi(\eta)\in\mathcal{P}(S)
\]
In addition, the so-called resolvent operators
\begin{equation}
P_{\eta}~:~f\in\mathcal{B}(S)\rightarrow P_{\eta}(f):=\sum_{n\geq0}\left[
M_{\eta}^{n}-\Phi(\eta)\right]  (f)\in\mathcal{B}(S)\label{defPeta}%
\end{equation}
are well defined absolutely convergent series that satisfy the Poisson
equation given by
\[
\left\{
\begin{array}
[c]{rcl}%
(M_{\eta}-Id)P_{\eta} & = & (\Phi(\eta)-Id)\\
\Phi(\eta)P_{\eta} & = & 0
\end{array}
\right.
\]

Resolvent operators are classical tools for the asymptotic analysis of time
inhomogeneous Markov chains. In our context the Markov chain interacts with a
flow a probability measures. To analyze the situation where this flow
converges to some limiting measure, it is convenient to study the regularity
properties of the resolvent operators $P_{\eta}$ as well as the ones of the
invariant measure mapping $\Phi(\eta)$ associated with $M_{\eta}$.

\begin{prop}
\label{propk} Under the regularity conditions (\ref{lipcnp1}) and
(\ref{lipcnp2}), we have
\begin{equation}
\sup_{\eta\in\mathcal{P}(S^{\prime})}\Vert P_{\eta}\Vert\leq p(n_{0}%
)\label{eq:first}%
\end{equation}
In addition, for any $f\in\mathcal{B}(S)$ and any $(\mu,\eta)\in
\mathcal{P}(S^{\prime})$ we have the following Lipschitz type inequalities
\begin{equation}
\left\vert \left[  \Phi(\eta)-\Phi(\mu)\right]  (f)\right\vert \leq
\int~\left\vert \left[  \eta-\mu\right]  (g)\right\vert \Upsilon_{\mu
}(f,dg)\label{o1}%
\end{equation}
and
\begin{equation}
\Vert\left[  P_{\eta}-P_{\mu}\right]  (f)\Vert\leq\int~\left\vert \left[
\eta-\mu\right]  (g)\right\vert \Upsilon_{\mu}^{\prime}(f,dg)\label{p1}%
\end{equation}
where $(\Upsilon_{\mu},\Upsilon_{\mu}^{\prime})$ is a pair of bounded integral
operators from $\mathcal{B}(S)$ into $\mathcal{B}(S^{\prime})$ indexed by the
set of measures $\mu\in\mathcal{P}(S^{\prime})$ such that
\[
\int~\Vert g\Vert~\Upsilon_{\mu}(f,dg)\leq p(n_{0})~\Lambda~\Vert f\Vert
\]
and
\[
\int~\Vert g\Vert~\Upsilon_{\mu}^{\prime}(f,dg)\leq p(n_{0})(1+p(n_{0}%
))~\Lambda~\Vert f\Vert
\]

\end{prop}

\noindent\mbox{\bf Proof:}\newline The first result (\ref{eq:first}) is proved
in~\cite{dd08}. For completeness, it is sketched here. We use the fact that
\[
P_{\eta}(f)(x)=\sum_{n\geq0}\int~\left[  M_{\eta}^{n}(f)(x)-M_{\eta}%
^{n}(f)(y)\right]  \Phi(\eta)(dy)
\]
to check that
\[
\Vert P_{\eta}(f)\Vert\leq\sum_{n\geq0}\mbox{\rm osc}(M_{\eta}^{n}(f))
\]
and
\[
\Vert P_{\eta}(f)\Vert\leq\left[  \sum_{n\geq0}\beta(M_{\eta}^{n})\right]
~\mbox{\rm osc}(f)\Rightarrow%
\begin{array}
[t]{rcl}%
\Vert P_{\eta}\Vert & \leq & 2~\sum_{p\geq1}~\sum_{r=0}^{n_{0}-1}\beta
(M_{\eta}^{(p-1)n_{0}+r})\\
& \leq & \frac{2n_{0}}{1-\beta(M_{\eta}^{n_{0}})}%
\end{array}
\]
The result (\ref{eq:first}) follows straightforwardly. The proof of (\ref{o1})
is based on the following decomposition
\[
\left[  \Phi(\eta)-\Phi(\mu)\right]  (f)=\left\{  \left[  \Phi(\eta)-\Phi
(\mu)\right]  M_{\mu}+\Phi(\eta)\left[  M_{\eta}-M_{\mu}\right]  \right\}
(f_{\mu})
\]
with $f_{\mu}:=(f-\Phi(\mu)(f))$. Under our regularity conditions on the
integral operators $M_{\mu}$, we find that
\begin{align*}
|\left[  \Phi(\eta)-\Phi(\mu)\right]  (f)| &  \leq|\left[  \Phi(\eta)-\Phi
(\mu)\right]  M_{\mu}(f_{\mu})|+\Vert\left[  M_{\eta}-M_{\mu}\right]  (f_{\mu
})\Vert\\
&  \leq|\left[  \Phi(\eta)-\Phi(\mu)\right]  M_{\mu}(f_{\mu})|+\int~\left\vert
\left[  \mu-\eta\right]  (g)\right\vert ~\Gamma_{\mu}(f_{\mu},dg)
\end{align*}
This recursion readily implies (\ref{o1}) with the integral operator given by
\[
\Upsilon_{\mu}(f,dg):=\sum_{n\geq0}\Gamma_{\mu}(M_{\mu}^{n}(f_{\mu}),dg)
\]
Finally we observe that
\[
\int~\Vert g\Vert~\Upsilon_{\mu}(f,dg)\leq\sum_{n\geq0}\int~\Vert
g\Vert~\Gamma_{\mu}(M_{\mu}^{n}(f_{\mu}),dg)\leq\Lambda~\sum_{n\geq0}\Vert
M_{\mu}^{n}(f_{\mu})\Vert
\]
Arguing as above, we conclude that
\[
\int~\Vert g\Vert~\Upsilon_{\mu}(f,dg)\leq\Lambda~\sum_{n\geq0}%
\mbox{\rm osc}(M_{\mu}^{n}(f))\leq p(n_{0})~\Lambda~\Vert f\Vert
\]
This ends the proof of (\ref{o1}). The proof of (\ref{p1}) follows the same
type of arguments. We observe that
\[
P_{\eta}-P_{\mu}=P_{\mu}(M_{\eta}-M_{\mu})P_{\eta}+\left[  \Phi(\mu)-\Phi
(\eta)\right]  P_{\eta}%
\]
To check this formula, we first use the fact that $M_{\mu}P_{\mu}=P_{\mu
}M_{\mu}$ to prove that
\[
P_{\mu}(M_{\mu}-Id)=(M_{\mu}-Id)P_{\mu}=(\Phi(\mu)-Id)
\]
This yields
\[
P_{\mu}(M_{\mu}-Id)P_{\eta}=(\Phi(\mu)-Id)P_{\eta}%
\]
Using the Poisson equation and the fact that $P_{\mu}(1)=0$ we also have the
decomposition
\[
P_{\mu}(M_{\eta}-Id)P_{\eta}=P_{\mu}(\Phi(\eta)-Id)=-P_{\mu}%
\]
Combining these two formulae, we conclude that
\[
P_{\mu}(M_{\eta}-M_{\mu})P_{\eta}=\left[  P_{\eta}-P_{\mu}\right]  -\left[
\Phi(\mu)-\Phi(\eta)\right]  P_{\eta}%
\]
This ends the proof of the decomposition given above. It is now easily checked
that
\begin{align*}
\Vert\left[  P_{\mu}-P_{\eta}\right]  (f)\Vert &  \leq\Vert(M_{\mu}-M_{\eta
})P_{\mu}(f)\Vert+\int~\left\vert \left[  \eta-\mu\right]  (g)\right\vert
\Upsilon_{\mu}(P_{\mu}(f),dg)\\
&  \leq\int~\left\vert \left[  \mu-\eta\right]  (g)\right\vert ~\left\{
\Gamma_{\mu}(P_{\mu}(f),dg)+\Upsilon_{\mu}(P_{\mu}(f),dg)\right\}
\end{align*}
The end of the proof follows the same type of arguments as before. This ends
the proof of the proposition. \hfill\hbox{\vrule height 5pt width 5pt depth
0pt}\medskip\newline

\subsection{Local fluctuations of weighted occupation measures}

\label{lln}

This section is concerned with the fluctuation analysis of the occupation
measures of the time inhomogeneous Markov chain introduced in (\ref{eqqM}). In
section~\ref{secfrf}, we shall use these results to analyze the fluctuations
of i-MCMC algorithms. The fluctuation analysis of this type of models is
related to the fluctuations of weighted occupation measures with respect to
some weight array type functions.

\begin{defi}
\label{defiW} We let $\mathcal{W}$ be the set of non negative and non
increasing weight array functions $w=(w_{n}(p))_{0\leq p\leq n,0\leq n}$,
satisfying the following conditions
\[
\exists m\geq1\quad\mbox{such that}\quad\sum_{n\geq0}w_{n}^{m}(0)<\infty
\]
with
\[
\forall\epsilon\in\lbrack0,1]\qquad\varpi(\epsilon):=\lim_{n\rightarrow\infty
}\sum_{0\leq p\leq\lfloor\epsilon n\rfloor}~w_{n}^{2}(p)<\infty
\]
and some scaling function $\varpi$ such that $\lim_{(\epsilon_{0},\epsilon
_{1})\rightarrow(0+,1-)}(\varpi(\epsilon_{0}),\varpi(\epsilon_{1}))=(0,1)$.
\end{defi}

We observe that the traditional and constant fluctuation rates sequences
$w_{n}(p)=1/\sqrt{n}$ belong to $\mathcal{W}$, with the identity function
$\varpi(\epsilon)=\epsilon$. In our setup it is necessary to introduce more
general sequences; see proof of theorem \ref{tcl} in Section \ref{secfrf}.

\begin{defi}
\label{defiWW} We associate to a given weight array function $w\in\mathcal{W}
$ the mapping
\[
W~:~\eta\in\mathcal{M}(S)^{\mathbb{N}}\mapsto W(\eta)=(W_{n}(\eta))_{n\geq
0}\in\mathcal{M}(S)^{\mathbb{N}}%
\]
defined for any flow of measures $\eta=(\eta_{n})_{n\geq0}\in\mathcal{P}(S)$,
and any $n\geq0$, by the weighted measures
\[
W_{n}(\eta)=\sum_{0\leq p\leq n}~w_{n}(p)~\eta_{p}%
\]

\end{defi}

The next proposition presents a pivotal decomposition formula of the weighted
occupation measures in terms of a martingale on fixed time horizon with a
negligible remainder bias term.

\begin{prop}
\label{pivotprop} \label{key} We consider the flow of random measures
$\zeta:=\left(  \zeta_{n}\right)  _{n\in\mathbb{N}}\in\mathcal{M}%
(S)^{\mathbb{N}}$ defined for any $n\geq0$ by the following formula
\[
\forall n\geq0\qquad\zeta_{n}=\left[  \delta_{X_{n}}-\Phi\left(  \eta
_{n-1}\right)  \right]
\]
For $n=0$, we use the convention $\Phi\left(  \eta_{-1}\right)  =\nu$ so that
$\zeta_{0}=\left[  \delta_{X_{0}}-\nu\right]  $. For any weight array function
$w\in\mathcal{W}$, the weighted measures $W_{n}(\zeta)$ satisfy the following
decomposition
\begin{equation}
W_{n}(\zeta)(f)=\sum_{0\leq p\leq n}w_{n}(p)~\Delta\mathbb{M}_{p+1}%
(f)+\mathbb{L}_{n}(f)\label{wzeta}%
\end{equation}
with the martingale increments
\begin{equation}
\Delta\mathbb{M}_{p+1}(f)=\mathbb{M}_{p+1}(f)-\mathbb{M}_{p}(f):=\left(
P_{\eta_{p-1}}(f)(X_{p+1})-M_{\eta_{p}}P_{\eta_{p-1}}(f)(X_{p})\right)
\label{marth}%
\end{equation}
and a remainder signed measure $\mathbb{L}_{n}$ such that
\[
\Vert\mathbb{L}_{n}\Vert\leq w_{n}(0)~(1+p(n_{0})\overline{\tau}%
(n))~(2+p(n_{0}))~\Lambda
\]

\end{prop}

\noindent\mbox{\bf Proof:}\newline We let $P_{\eta_{n-1}}$ be the integral
operator solution of the Poisson equation associated with the Markov
transition $M_{\eta_{n-1}}$ with an invariant measure $\Phi\left(  \eta
_{n-1}\right)  $. By construction, we have
\[
\zeta_{n}(f)=\left[  f(X_{n})-\Phi\left(  \eta_{n-1}\right)  (f)\right]
=P_{\eta_{n-1}}(f)(X_{n})-M_{\eta_{n-1}}(P_{\eta_{n-1}}(f))(X_{n})
\]
For $n=0$, we use the convention $P_{\eta_{-1}}=Id$ and $M_{\eta_{-1}}=\nu$.
The proof of (\ref{wzeta}) is based on the following decomposition
\[
\zeta_{n}(f)=A_{n}(f)+B_{n}(f)+C_{n}(f)+\Delta\mathbb{M}_{n+1}(f)
\]
with the random processes $A_{n}(f)$, $B_{n}(f)$ and $C_{n}(f)$ defined below
\begin{align*}
A_{n}(f):= &  \left[  P_{\eta_{n}}-P_{\eta_{n-1}}\right]  (f)(X_{n+1})\\
B_{n}(f):= &  \left[  P_{\eta_{n-1}}(f)(X_{n})-P_{\eta_{n}}(f)(X_{n+1})\right]
\\
C_{n}(f):= &  \left[  M_{\eta_{n}}-M_{\eta_{n-1}}\right]  P_{\eta_{n-1}%
}(f)(X_{n})
\end{align*}
Using the Lipschitz inequality (\ref{p1}) presented in proposition~\ref{propk}%
, we prove that
\begin{align*}
\left\vert A_{n}(f)\right\vert  &  \leq\Vert\left[  P_{\eta_{n}}-P_{\eta
_{n-1}}\right]  (f)\Vert\\
&  \leq\int~\left\vert \left[  \eta_{n}-\eta_{n-1}\right]  (g)\right\vert
\Upsilon_{\eta_{n}}^{\prime}(f,dg)\leq\tau(n)~p(n_{0})(1+p(n_{0}%
))~\Lambda~\Vert f\Vert
\end{align*}
In addition, using the Lipschitz regularity condition (\ref{lipcnp2}), we also
obtain%
\begin{align*}
\left\vert C_{n}(f)\right\vert  &  \leq\Vert\left[  M_{\eta_{n}}-M_{\eta
_{n-1}}\right]  (P_{\eta_{n-1}}(f))\Vert\\
&  \leq\Vert\int~\left\vert \left[  \eta_{n}-\eta_{n-1}\right]  (g)\right\vert
~\Gamma_{\eta_{n}}(P_{\eta_{n-1}}(f),dg)\\
&  \leq\tau(n)~\Lambda~\Vert P_{\eta_{n-1}}\Vert~\Vert f\Vert\leq
\tau(n)~\Lambda~p(n_{0})~\Vert f\Vert
\end{align*}
By definition of the weighted measure $W_{n}(\zeta)$, we have
\begin{align}
W_{n}(\zeta):= &  \sum_{0\leq p\leq n}w_{n}(p)\zeta_{p}(f)\nonumber\\
&  =\sum_{0\leq p\leq n}w_{n}(p)~\Delta\mathbb{M}_{p+1}(f)+\sum_{0\leq p\leq
n}w_{n}(p)(A_{p}(f)+B_{p}(f)+C_{p}(f))\label{refp1}%
\end{align}
From previous calculations, we have
\[
\left\vert \sum_{0\leq p\leq n}w_{n}(p)(A_{p}(f)+C_{p}(f))\right\vert \leq
w_{n}(0)~\overline{\tau}(n)~\Lambda~p(n_{0})(2+p(n_{0}))~\Vert f\Vert
\]
Finally, we use the following decomposition
\begin{align*}
\sum_{0\leq p\leq n}w_{n}(p)B_{p}(f) &  =\sum_{0\leq p\leq n}\left[
w_{n}(p)~P_{\eta_{p-1}}(f)(X_{p})-w_{n}(p+1)P_{\eta_{p}}(f)(X_{p+1})\right] \\
&  \hskip2cm+\sum_{0\leq p\leq n}\left[  w_{n}(p+1)-w_{n}(p)\right]
~P_{\eta_{p}}(f)(X_{p+1})
\end{align*}
with the convention $w_{n}(n+1)=0$. This implies that
\begin{align*}
\left\vert \sum_{0\leq p\leq n}w_{n}(p)B_{p}(f)\right\vert  &  \leq
2~w_{n}(0)~\Vert f\Vert+p(n_{0})~\Vert f\Vert~\sum_{0\leq p\leq n}\left[
w_{n}(p)-w_{n}(p+1)\right] \\
&  =(2+p(n_{0}))~\Vert f\Vert~w_{n}(0)
\end{align*}
The end of the proof is now a direct consequence of formula (\ref{refp1}).
\hfill\hbox{\vrule height 5pt width 5pt depth 0pt}\medskip\newline

Now, we are in position to state and to prove the main result of this section.

\begin{theo}
\label{theointernon}Assume there exist a measure $\eta$ and some $m\geq1$ such
that
\[
\forall f\in\mathcal{B}_{1}(S^{\prime})\qquad\mathbb{E}(|\eta_{n}%
(f)-\eta(f)|^{m})\leq\epsilon_{m}(n)\quad\quad\mbox{\rm with}\quad\sum
_{n\geq0}\epsilon_{m}(n)<\infty
\]
We let $V_{n}:=W_{n}(\zeta)$ be the sequence of random fields on
$\mathcal{B}(S)$ associated with a given weight array function $w\in
\mathcal{W}$ and defined in (\ref{wzeta}). We suppose that $w\in\mathcal{W}$
is chosen so that $w_{n}(0)\overline{\tau}(n)$ tends to $0$ as $n\rightarrow
\infty$. In this situation, $V_{n}$ converges in law as $n\rightarrow\infty$
to a Gaussian random field $V$ on $\mathcal{B}(S)$ such that
\[
\forall(f,g)\in\mathcal{B}(S)^{2}\qquad\mathbb{E}(V(f)V(g))=\Phi(\eta)\left[
C_{\eta}(f,g)\right]
\]
with the local covariance function
\[
C_{\eta}(f,g):=M_{\eta}\left[  \left(  P_{\eta}(f)-M_{\eta}P_{\eta}(f)\right)
\left(  P_{\eta}(g)-M_{\eta}P_{\eta}(g)\right)  \right]
\]

\end{theo}

\noindent\mbox{\bf Proof:}\newline Using proposition~\ref{key}, it is clearly
sufficient to prove that the random fields
\begin{equation}
W_{n}^{\prime}(\zeta):=\sum_{0\leq p\leq n}w_{n}(p)~\Delta\mathbb{M}%
_{p+1}\label{randomfieldprime}%
\end{equation}
converge in law to the Gaussian random field $V$ as $n\rightarrow\infty$. To
use the Lindeberg CLT for triangular arrays of $\mathbb{R}^{d}$-valued random
variables, we let $f=(f^{i})_{1\leq i\leq d}\in\mathcal{B}(S)^{d}$ be a
collection of $d$-valued functions and we consider the $\mathbb{R}^{d}$-valued
random variables $W_{n}^{\prime}(\zeta)(f)=(W_{n}^{\prime}(\zeta
)(f^{i}))_{1\leq i\leq d}$. We further denote by $\mathcal{F}_{p}$ the
$\sigma$-field generated by the random variables $X_{q}$ for any $q\leq p$. By
construction, for any functions $f$ and $g\in\mathcal{B}(S)$ and for every
$0\leq p\leq n$ we find that
\begin{align*}
\mathbb{E}(w_{n}(p)~\Delta\mathbb{M}_{p+1}(f)~|~\mathcal{F}_{p}) &  =0\\
\mathbb{E}(w_{n}(p)^{2}~\Delta\mathbb{M}_{p+1}(f)\Delta\mathbb{M}%
_{p+1}(g)~|~\mathcal{F}_{p}) &  =w_{n}(p)^{2}~C_{p}^{\prime}(f,g)(X_{p})
\end{align*}
with the local covariance function
\[%
\begin{array}
[c]{l}%
C_{p}^{\prime}(f,g):=M_{\eta_{p}}\left[  \left(  P_{\eta_{p-1}}(f)-M_{\eta
_{p}}P_{\eta_{p-1}}(f)\right)  \left(  P_{\eta_{p-1}}(g)-M_{\eta_{p}}%
P_{\eta_{p-1}}(g)\right)  \right]
\end{array}
\]
Using proposition~\ref{propk}, after some tedious but elementary calculations
we find that
\begin{align*}
\Vert C_{p}^{\prime}(f,g)-C_{\eta}(f,g)\Vert &  \leq c(\eta)~\left\{
\int\left\vert \left[  \eta_{p-1}-\eta\right]  (h)\right\vert \Upsilon_{\eta
}^{1}((f,g),dh)\right. \\
&  \hskip3cm+\left.  \int\left\vert \left[  \eta_{p}-\eta\right]
(h)\right\vert \Upsilon_{\eta}^{2}((f,g),dh)\right\}
\end{align*}
with a pair of bounded integral operator $\Upsilon_{\eta}^{i}$, $i=1,2$, from
$\mathcal{B}(S)^{2}$ into $\mathcal{B}(S^{\prime})$ such that
\[
\int~\Vert h\Vert~\Upsilon_{\eta}^{i}((f,g),dh)\leq c(\eta)~\Vert f\Vert\Vert
g\Vert
\]
In the above displayed formula, $c(\eta)<\infty$ stands for a finite constant
whose value only depends on the measure $\eta$. Under our assumptions, the
following almost sure convergence result readily follows
\begin{equation}
\lim_{p\rightarrow\infty}\Vert C_{p}^{\prime}(f,g)-C_{\eta}(f,g)\Vert
=0\label{kk1}%
\end{equation}
On the other hand, using proposition~\ref{pivotprop}, for any function
$h\in\mathcal{B}(S)$ we have the decomposition
\[
\sum_{0\leq p\leq n}w_{n}^{2}(p)~(h(X_{p})-\Phi(\eta_{p-1})(h))=\sum_{0\leq
p\leq n}w_{n}^{2}(p)~\Delta\mathbb{M}_{p+1}(h)+\mathbb{L}_{n}^{\prime}(h)
\]
with the martingale increments $\Delta\mathbb{M}_{p+1}(h)$ given in
(\ref{marth}), and a remainder signed measure $\mathbb{L}_{n}^{\prime}$ such
that
\[
\Vert\mathbb{L}_{n}^{\prime}\Vert\leq w_{n}^{2}(0)~(1+p(n_{0})\overline{\tau
}(n))~(2+p(n_{0}))~\Lambda
\]
Using a Burkholder-Davis-Gundy type inequality for martingales, we find that
for any $m\geq1$
\[
\mathbb{E}\left(  \left\vert \sum_{0\leq p\leq n}w_{n}^{2}(p)~\Delta
\mathbb{M}_{p+1}(h)\right\vert ^{m}\right)  ^{\frac{1}{m}}\leq a(m)~\left(
\sum_{0\leq p\leq n}~w_{n}(p)^{4}\right)  ^{\frac{1}{2}}~\mbox{\rm osc}(h)
\]
for some constants $a(m)$ whose values only depend on the parameter $m$. Thus,
we find that
\[
\mathbb{E}\left(  \left\vert \sum_{0\leq p\leq n}w_{n}^{2}(p)~\Delta
\mathbb{M}_{p+1}(h)\right\vert ^{m}\right)  ^{\frac{1}{m}}\leq a(m)~w_{n}%
(0)~\left(  \sum_{0\leq p\leq n}~w_{n}(p)^{2}\right)  ^{\frac{1}{2}%
}~\mbox{\rm osc}(h)
\]
Under our assumptions on the weight functions $w$, if we take $h=C_{\eta
}(f,g)$ then by (\ref{kk1}) we obtain the following almost sure convergence
result
\[
\lim_{n\rightarrow\infty}\sum_{p=0}^{n}w_{n}(p)^{2}~C_{\eta}(f,g)(X_{p}%
)=\lim_{n\rightarrow\infty}\sum_{p=0}^{n}w_{n}(p)^{2}~\Phi(\eta_{p-1}%
)(C_{\eta}(f,g))
\]
We now combine the regularity property (\ref{o1}) with the generalized
Minkowski inequality to prove that for any function $h\in\mathcal{B}(S)$ and
any $m\geq1$
\[%
\begin{array}
[c]{l}%
\mathbb{E}\left(  \left\vert \sum_{p=0}^{n}w_{n}(p)^{2}~[\Phi(\eta
_{p-1})(h)-\Phi(\eta)(h)]\right\vert ^{m}\right)  ^{\frac{1}{m}}\\
\\
\leq\displaystyle\sum_{p=0}^{n}w_{n}(p)^{2}~\int~\mathbb{E}\left(  \left\vert
\left[  \eta_{p-1}-\eta\right]  (g)\right\vert ^{m}\right)  ^{\frac{1}{m}%
}\Upsilon_{\mu}(h,dg)\\
\\
\leq(p(n_{0})\Lambda\Vert h\Vert)~\displaystyle\left(  w_{n}(0)^{2}+\sum
_{p=1}^{n}w_{n}(p)^{2}\epsilon_{m}(p-1)\right)
\end{array}
\]
This readily implies that
\[%
\begin{array}
[c]{l}%
\mathbb{E}\left(  \left\vert \sum_{p=0}^{n}w_{n}(p)^{2}~[\Phi(\eta
_{p-1})(h)-\Phi(\eta)(h)]\right\vert ^{m}\right)  ^{\frac{1}{m}}\\
\\
\leq w_{n}(0)^{2}~(p(n_{0})\Lambda\Vert h\Vert)~\left(  1+\displaystyle\sum
_{p\geq0}\epsilon_{m}(p)\right)
\end{array}
\]
If we choose $h=C_{\eta}(f,g)$, this yields the following almost sure
convergence result
\[
=\lim_{n\rightarrow\infty}\sum_{p=0}^{n}w_{n}(p)^{2}~\Phi(\eta_{p-1})(C_{\eta
}(f,g))=\Phi(\eta)(C_{\eta}(f,g))
\]
To summarize, we have proved the following series of almost sure convergence
results
\begin{align*}
\lim_{n\rightarrow\infty}\sum_{p=0}^{n}w_{n}(p)^{2}~C_{p}^{\prime}(f,g)(X_{p})
&  =\lim_{n\rightarrow\infty}\sum_{p=0}^{n}w_{n}(p)^{2}~C_{\eta}(f,g)(X_{p})\\
&  =\lim_{n\rightarrow\infty}\sum_{p=0}^{n}w_{n}(p)^{2}~\Phi(\eta
_{p-1})(C_{\eta}(f,g))\\
&  =\Phi(\eta)(C_{\eta}(f,g))
\end{align*}
Therefore, we also have the almost sure convergence result
\[
\lim_{n\rightarrow\infty}\sum_{p=0}^{n}w_{n}(p)^{2}~\mathbb{E}(\Delta
\mathbb{M}_{p+1}(f)\Delta\mathbb{M}_{p+1}(g)~|~\mathcal{F}_{p})=\Phi
(\eta)(C_{\eta}(f,g))
\]
Since we have $\vee_{0\leq p\leq n}w_{n}(p)=w_{n}(0)\rightarrow0$, as
$n\rightarrow\infty$, the Lindeberg condition is satisfied and we conclude
that the sequence of random fields $W_{n}^{\prime}(\zeta)$ defined in
(\ref{randomfieldprime}) converges in law to the Gaussian random field $V$ as
$n\rightarrow\infty$. This ends the proof of the theorem. \hfill\hbox{\vrule
height 5pt width 5pt depth 0pt}\medskip\newline

We end this section with an alternative and simpler representation of the
covariance function of the random field $V$ presented in
theorem~\ref{theointernon}.\ We have
\[
C_{\eta}(f,f)(x)=\int~M_{\eta}(x,dy)~\left[  P_{\eta}(f)(y)-M_{\eta}(P_{\eta
}(f))(x)\right]  ^{2}%
\]
Using the decomposition
\begin{align*}
P_{\eta}(f)(y)-M_{\eta}(P_{\eta}(f))(x)  &  =\left[  P_{\eta}(f)(y)-P_{\eta
}(f)(x)\right]  +\left[  P_{\eta}(f)(x)-M_{\eta}(P_{\eta}(f))(x)\right] \\
&  =\left[  P_{\eta}(f)(y)-P_{\eta}(f)(x)\right]  +\left[  f(x)-\Phi
(\eta)(f)\right]
\end{align*}
and the fact that
\begin{align*}
\int~M_{\eta}(x,dy)~\left[  P_{\eta}(f)(y)-P_{\eta}(f)(x)\right]   &  =\left[
M_{\eta}(P_{\eta}(f))(x)-P_{\eta}(f)(x)\right] \\
&  =-\left[  f(x)-\Phi(\eta)(f)\right]
\end{align*}
we prove the formula
\[
C_{\eta}(f,f)(x)=\int~M_{\eta}(x,dy)~\left[  P_{\eta}(f)(y)-P_{\eta
}(f)(x)\right]  ^{2}-\left[  f(x)-\Phi(\eta)(f)\right]  ^{2}%
\]
On the other hand, recalling that $\Phi(\eta)=\Phi(\eta)M_{\eta}$ and using
again the Poisson equation we also have
\[
\int~\Phi(\eta)(dx)M_{\eta}(x,dy)~\left[  P_{\eta}(f)(y)-P_{\eta
}(f)(x)\right]  ^{2}=2~\Phi(\eta)\left[  P_{\eta}(f)~(f-\Phi(\eta)(f))\right]
\]
and
\[%
\begin{array}
[c]{l}%
2~\Phi(\eta)\left[  P_{\eta}(f)~(f-\Phi(\eta)(f))\right] \\
\\
=2~\Phi(\eta)\left[  (f-\Phi(\eta)(f))^{2}\right]  +2\sum_{n\geq1}~\Phi
(\eta)\left[  M_{\eta}^{n}(f-\Phi(\eta)(f))~(f-\Phi(\eta)(f))\right]
\end{array}
\]
Hence we have proved the following proposition.

\begin{prop}
The limiting covariance function presented in theorem~\ref{theointernon} is
alternatively defined for any function $f\in\mathcal{B }(S)$ by the following
formula
\[%
\begin{array}
[c]{l}%
\Phi(\eta) \left[  C_{\eta}(f,f)\right] \\
\\
= \Phi(\eta)\left[  (f- \Phi(\eta)(f))^{2}\right]  + 2\sum_{n\geq1}~ \Phi
(\eta)\left[  (f- \Phi(\eta)(f))~M^{n}_{\eta}(f- \Phi(\eta)(f))\right]
\end{array}
\]

\end{prop}

\section{A fluctuation theorem for local interaction fields}

\label{secfrf}

\subsection{Introduction}

\label{introf} This section presents the fluctuation analysis of a class of
weighted random fields associated to i-MCMC algorithms. Following the local
fluctuation analysis for time inhomogeneous Markov chains presented in
section~\ref{lln}, we introduce the following weighted random fields.

\begin{defi}
We consider the flow of random measures
\[
\forall l\geq0\quad\forall n\geq0\qquad\delta_{n}^{(l)}:=\left[  \delta
_{X_{n}^{(l)}}-\Phi^{\left(  l\right)  }\left(  \eta_{n-1}^{(l-1)}\right)
\right]
\]
For $n=0$, we use the convention $\Phi^{\left(  l\right)  }\left(  \eta
_{-1}\right)  =\nu^{\left(  l\right)  }$ so that $\delta_{0}^{(l)}=\left[
\delta_{X_{0}^{(l)}}-\nu^{\left(  l\right)  }\right]  $. We associate to a
sequence of weight array functions $(w^{(l)})_{l\geq0}\in\mathcal{W}%
^{\mathbb{N}}$ the following flow of random fields $(W_{n}^{(l)}(\delta
^{(l)}))_{l\geq0}$ on the sets of functions $(\mathcal{B}(S^{(l)}))_{l\geq0}$
\[
\forall l\geq0\quad\forall n\geq0\qquad W_{n}^{(l)}(\delta^{(l)}):=\sum_{0\leq
p\leq n}w_{n}^{(l)}(p)~\delta_{p}^{(l)}%
\]

\end{defi}

We observe that the regularity conditions (\ref{lipcna}) and (\ref{lipcnb})
ensure that the collection of Markov operators $M_{\eta}^{(l)}$ and their
invariant measures $\Phi^{\left(  l\right)  }(\eta)$ satisfy the regularity
conditions (\ref{lipcnp1}) and (\ref{lipcnp2}) introduced in
section~\ref{dmmc}. Also observe that the i-MCMC chain $X^{(l+1)}$ is a time
inhomogeneous model of the form (\ref{eqqM}) with a collection of transitions
$M_{\eta_{n}^{(l)}}^{(l+1)}$ that depend on the flow of occupation measures
$\eta_{n}^{(l)}$ associated with the i-MCMC chain at level $l$. In this
scenario, condition (\ref{lipcnp3}) is satisfied with
\[
\forall n\geq0\qquad\Vert\eta_{n}^{(l)}-\eta_{n-1}^{(l)}\Vert\leq\tau
^{(l)}(n):=\frac{2}{n+1}%
\]
and we have
\[
\overline{\tau}^{(l)}(n):=\sum_{0\leq p\leq n}\tau^{(l)}(p)=2\sum_{0\leq p\leq
n}\frac{1}{p+1}\leq2(1+\log{(n+1)})
\]
Finally, we recall that for any $m\geq1$ we have
\begin{equation}
\forall l\geq0\quad\forall f\in\mathcal{B}_{1}(S^{(l)})\qquad\mathbb{E}%
(|\eta_{n}^{(l)}(f)-\pi^{(l)}(f)|^{m})^{\frac{1}{m}}\leq b(m)~c(l)~\frac
{1}{\sqrt{n+1}}\label{dlp}%
\end{equation}
for a collection of finite constants $b(m)$ whose values only depend on the
parameter $m$ (see for instance~\cite{dd08}). Using theorem~\ref{theointernon}%
, we can prove that the random fields
\begin{equation}
V_{n}^{(l)}:=W_{n}^{(l)}(\delta^{(l)})\label{randF}%
\end{equation}
associated with a given weight array function $w^{(l)}\in\mathcal{W}$
converges in law to a Gaussian random field $V^{(l)}$ as $n\rightarrow\infty$
such that
\begin{equation}
\forall(f,g)\in\mathcal{B}(S^{(l)})^{2}\qquad\mathbb{E}(V^{(l)}(f)V^{(l)}%
(g))=\pi^{(l)}\left[  C^{(l)}(f,g)\right] \label{defcov}%
\end{equation}
Here the covariance functions $C^{(l)}(f,g)$ are defined in terms of the
resolvent operator $P_{\pi^{(l-1)}}^{(l)}$ associated to the Markov transition
$M_{\pi^{(l-1)}}^{(l)}$ and the fixed point measure $\Phi^{\left(  l\right)
}(\pi^{(l-1)})=\pi^{(l)}$ with the following formula
\[%
\begin{array}
[c]{l}%
C^{(l)}(f,g)\\
\\
:=M_{\pi^{(l-1)}}^{(l)}\left[  \left(  P_{\pi^{(l-1)}}^{(l)}(f)-M_{\pi
^{(l-1)}}^{(l)}P_{\pi^{(l-1)}}^{(l)}(f)\right)  \left(  P_{\pi^{(l-1)}}%
^{(l)}(g)-M_{\pi^{(l-1)}}^{(l)}P_{\pi^{(l-1)}}^{(l)}(g)\right)  \right]
\end{array}
\]

The main objective of this section is to prove the following theorem.

\begin{theo}
\label{thf} We consider a collection of weight array functions $(w^{(l)}%
)_{l\geq0}\in\mathcal{W}^{\mathbb{N}}$. In this situation, the corresponding
flow of weighted random fields $(V_{n}^{(l)})_{l\geq0}$ defined in
(\ref{randF}), converges in law, as $n$ tends to infinity and in the sense of
finite dimensional distributions, to a sequence of independent and centered
Gaussian fields $\left(  V^{(l)}\right)  _{l\geq0}$ with covariance functions
defined in (\ref{defcov}).
\end{theo}

Using this result, the proof of the multivariate centra limit
theorem~\ref{tcl} follows exactly the same arguments as the ones we used in
the proof of theorem 2.1 in~\cite[section 6]{bdd08}.

\textbf{Proof of theorem~\ref{tcl} :}\newline We let $\mathbb{S}%
^{k}:=\mathbb{S}\mathbb{S}^{k-1}$ be the $k$-th iterate of the mapping
$\mathbb{S}:\eta\in\mathcal{M}(S^{(l)})^{\mathbb{N}}\mapsto\mathbb{S}%
(\eta)=(\mathbb{S}_{n}(\eta))_{n\geq0}\in\mathcal{M}(S^{(l)})^{\mathbb{N}}$
defined for any $\eta=(\eta_{n})_{n\geq0}\in\mathcal{M}(S^{(l)})^{\mathbb{N}}$
by
\[
\forall n\geq0\qquad\mathbb{S}_{n}(\eta)=\frac{1}{n+1}\sum_{0\leq p\leq
n}~\eta_{p}%
\]
We observe that the time averaged semigroup $\mathbb{S}^{k}$ can be rewritten
in terms of the following weighted summations
\[
\mathbb{S}_{n}^{k}(\eta)=\frac{1}{n+1}\sum_{0\leq p\leq n}~s_{n}^{(k)}%
(p)~\eta_{p}%
\]
with the weight array functions $s_{n}^{(k)}:=(s_{n}^{(k)}(p))_{0\leq p\leq
n}$ defined by
\[
\forall k\geq1\quad\forall0\leq p\leq n\qquad s_{n}^{(k+1)}(p)=\sum_{p\leq
q\leq n}~\frac{1}{(q+1)}~s_{n}^{(k)}(q)\quad\mbox{\rm and}\quad s_{n}%
^{(1)}(p):=1
\]
We also know from proposition 6.1 in~\cite{bdd08} that
\[
\lim_{n\rightarrow\infty}\frac{1}{n}\sum_{0\leq q\leq n}s_{n}^{(k+1)}%
(q)^{2}={(2k)!}/{k!^{2}}%
\]
and, for any $k\geq1$, the weight array functions $w^{(k)}$ defined by
\[
\forall n\geq0\quad\forall0\leq p\leq n\qquad w_{n}^{(k)}(p):={s_{n}^{(k)}%
(p)}/\sqrt{\sum_{0\leq q\leq n}s_{n}^{(k)}(q)^{2}}%
\]
belong to the set $\mathcal{W}$ introduced in definition~\ref{defiW}.

Moreover, using proposition 5.2 in~\cite{bdd08}, we have the following
multilevel expansion
\begin{equation}
\eta_{n}^{(k)}-\pi^{(k)}=\sum_{0\leq l\leq k}\mathbb{S}_{n}^{(l+1)}%
(\delta^{(k-l)})~D_{(k-l)+1,k}+\Xi_{n}^{(k)}\label{laform}%
\end{equation}
where $\Xi^{(k)}=(\Xi_{n}^{(k)})_{n\geq0}$ is a flow of signed random measures
such that
\[
\forall m\geq1\qquad\sup_{f\in\mathcal{B}_{1}(S^{(k)})}\mathbb{E}(|\Xi
_{n}^{(k)}(f)|^{m})^{\frac{1}{m}}\leq b(m)~c(k)~{\left(  \log{(n+1)}\right)
^{k}}/{(n+1)}%
\]
Here $b(m)$ stands for some constant whose value only depends on the parameter
$m$. This multilevel expansion implies that
\[%
\begin{array}
[c]{l}%
\sqrt{(n+1)}\left[  \eta_{n}^{(k)}-\pi^{(k)}\right] \\
\\
=\sum_{0\leq l\leq k}\sqrt{\frac{1}{n+1}\sum_{0\leq q\leq n}s_{n}%
^{(l+1)}(q)^{2}}~W_{n}^{(k-l)}(\delta^{(k-l)})~D_{(k-l)+1,k}+\overline{\Xi
}_{n}^{(k)}%
\end{array}
\]
with the weighted distribution flow mappings $W^{(k-l)}$ associated to the
weight functions $w^{(l+1)}$ and a remainder signed measure $\overline{\Xi
}_{n}^{(k)}$ such that
\[
\sup_{f\in\mathcal{B}_{1}(S^{(k)})}\mathbb{E}(|\overline{\Xi}_{n}%
^{(k)}(f)|)\leq c(k)~{\left(  \log{(n+1)}\right)  ^{k}}/\sqrt{(n+1)}%
\]
The proof of theorem~\ref{tcl} is now a direct consequence of
theorem~\ref{thf}. \hfill\hbox{\vrule height 5pt
width 5pt depth 0pt}\medskip\newline

\subsection{A martingale limit theorem}

This section is mainly concerned with the proof of theorem~\ref{thf}. We
following the same lines of arguments as the ones used in section~\ref{lln}
devoted to the fluctuations of weighted occupation measures associated with
time inhomogeneous Markov chains.

First, we introduce a few notation. For any $k\geq0$ and any $\mu
\in\mathcal{P}(S^{(k-1)})$, let $P_{\mu}^{(k)}$ be the resolvent operator
associated to the Markov transition $M_{\mu}^{(k)}$ and its invariant measure
$\Phi^{\left(  k\right)  }(\mu)\in\mathcal{P}(S^{(k)})$. We also set
\[
p^{(k)}(n_{k}):={2n_{k}}/{(1-m_{k}(n_{k}))}%
\]
with the pair of parameters $(n_{k},m_{k})$ defined in (\ref{lipcna}).

Using proposition~\ref{key}, we find that the weighted measures $W_{n}%
^{(k)}(\delta^{(k)})$ satisfy the following decomposition
\[
W_{n}^{(k)}(\delta^{(k)})(f)=\sum_{0\leq p\leq n}w_{n}^{(k)}(p)~\Delta
\mathbb{M}_{p+1}^{(k)}(f)+\mathbb{L}_{n}^{(k)}(f)
\]
for any $f\in\mathcal{B}(S^{(k)})$ with the martingale increments
\[
\Delta\mathbb{M}_{p+1}^{(k)}(f)=\mathbb{M}_{p+1}^{(k)}(f)-\mathbb{M}_{p}%
^{(k)}(f)=\left(  P_{\eta_{p-1}^{(k-1)}}^{(k)}(f)(X_{p+1}^{(k)})-M_{\eta
_{p}^{(k-1)}}^{(k)}P_{\eta_{p-1}^{(k-1)}}^{(k)}(f)(X_{p}^{(k)})\right)
\]
and the remainder signed measure $\mathbb{L}_{n}^{(k)}$ which are such that
\[
\Vert\mathbb{L}_{n}^{(k)}\Vert\leq\left\{  w_{n}^{(k)}(0)~\left(
1+2p^{(k)}(n_{k})(1+\log{(n+1)})\right)  ~(2+p^{(k)}(n_{k}))~\Lambda\right\}
\longrightarrow_{n\rightarrow\infty}0
\]
We consider a sequence of functions $f=(f^{i})_{1\leq i\leq d}$, with $d\geq
1$, and $f^{i}=(f_{k}^{i})_{k\geq0}\in\prod_{k\geq0}\mathcal{B}(S^{(k)})$, and
we let $\mathcal{W}^{(n)}(f)=(\mathcal{W}^{(n)}(f^{i}))_{1\leq i\leq d} $ be
the $\mathbb{R}^{d}$-valued and $\mathcal{F}^{(n)}$-adapted process defined
for any $l\geq0$ and any $1\leq i\leq d$ by
\[
\mathcal{W}_{l}^{(n)}(f^{i}):=\sum_{0\leq k\leq l}W_{n}^{(k)}(\delta
^{(k)})(f_{k}^{i})
\]
From the previous discussion, we find that
\[
\mathcal{W}_{l}^{(n)}(f^{i})=\mathcal{M}_{l}^{(n)}(f^{i})+\mathcal{L}%
_{l}^{(n)}(f^{i})
\]
with the $\mathcal{F}^{(n)}$-martingale $\mathcal{M}_{l}^{(n)}(f^{i})$ given
below
\begin{equation}
\mathcal{M}_{l}^{(n)}(f^{i}):=\sum_{0\leq k\leq l}\Delta\mathcal{M}_{k}%
^{(n)}(f^{i})\quad\mbox{\rm with}\quad\Delta\mathcal{M}_{k}^{(n)}(f^{i}%
):=\sum_{0\leq p\leq n}w_{n}^{(k)}(p)~\Delta\mathbb{M}_{p+1}^{(k)}%
(f)\label{martdef}%
\end{equation}
and the remainder bias type measure $\mathcal{L}_{l}^{(n)}=\sum_{0\leq k\leq
l}\mathbb{L}_{n}^{(k)}$, such that
\[
\lim_{n\rightarrow\infty}\Vert\mathcal{L}_{l}^{(n)}\Vert=0
\]

Theorem~\ref{thf} is now a direct consequence of the following proposition
(see for instance the arguments used in section 4.2 in~\cite{bdd08}).

\begin{prop}
\label{theointer} The sequence of martingales $\mathcal{M}_{l}^{(n)}(f)$
defined in (\ref{martdef}) converges in law as $n\rightarrow\infty$ to an
$\mathbb{R}^{d}$-valued Gaussian martingale $\mathcal{M}_{l}(f)=(\mathcal{M}%
_{l}(f^{i}))_{1\leq i\leq d}$ such that for any $l\geq0$ and any pair of
indexes $1\leq i,j\leq d$
\[
\langle\mathcal{M}(f^{i}),\mathcal{M}(f^{j})\rangle_{l}=\sum_{0\leq k\leq
l}\pi^{(k)}\left[  C^{(k)}(f^{i},f^{j})\right]
\]
with the local covariance functions $\pi^{(k)}\left[  C^{(k)}(f^{i}%
,f^{j})\right]  $ defined in (\ref{defcov}).
\end{prop}

\noindent\mbox{\bf Proof:}\newline The proof of the proposition is along the
same lines as the proof of theorem~\ref{theointernon}. First, we consider the
decomposition
\[
\mathcal{M}_{l}^{(n)}(f^{i})=\sum_{i=0}^{l(n+1)+n}~\mathcal{V}_{i}^{(n)}(f)
\]
where for every $0\leq i\leq l(n+1)+n$, with $i=k(n+1)+p$ for some $0\leq
k\leq l$, and $0\leq p\leq n$
\[
\mathcal{V}_{i}^{(n)}(f):=w_{n}^{(k)}(p)~\Delta\mathbb{M}_{p+1}^{(k)}(f_{k})
\]
We further denote by $\mathcal{G}_{i}^{(n)}$ the $\sigma$-field generated by
the pair of random variables $(X_{p}^{(k)},X_{p+1}^{(k)})$ for any pair of
parameters $(k,p)$ such that $k(n+1)+p\leq i$. By construction, for any flow
of functions $f=(f_{l})_{l\geq0}$ and $g=(g_{l})_{l\geq0}\in\prod_{l\geq
0}\mathcal{B}(S^{(l)})$ and for every $0\leq i\leq l(n+1)+n$, with
$i=k(n+1)+p$ for some $0\leq k\leq l$, and $0\leq p\leq n$, we find that
\begin{align*}
\mathbb{E}(\mathcal{V}_{i}^{(n)}(f)~|~\mathcal{G}_{i-1}^{(n)}) &  =0\\
\mathbb{E}(\mathcal{V}_{i}^{(n)}(f)\mathcal{V}_{i}^{(n)}(g)~|~\mathcal{G}%
_{i-1}^{(n)}) &  =w_{n}^{(k)}(p)^{2}~C_{p}^{(k)}(f,g)(X_{p}^{(k)})
\end{align*}
with the local covariance function
\[%
\begin{array}
[c]{l}%
C_{p}^{(k)}(f,g):=M_{\eta_{p}^{(k-1)}}^{(k)}\left[  \left(  P_{\eta
_{p-1}^{(k-1)}}^{(k)}(f_{k})-M_{\eta_{p}^{(k-1)}}^{(k)}P_{\eta_{p-1}^{(k-1)}%
}^{(k)}(f_{k})\right)  \right. \\
\\
\left.  \hskip6cm\times\left(  P_{\eta_{p-1}^{(k-1)}}^{(k)}(g_{k})-M_{\eta
_{p}^{(k-1)}}^{(k)}P_{\eta_{p-1}^{(k-1)}}^{(k)}(g_{k})\right)  \right]
\end{array}
\]
Under our Lipschitz regularity conditions (\ref{lipcna}) and (\ref{lipcnb}),
proposition~\ref{propk} applies to the mappings $\Phi^{\left(  k\right)  }$
and the resolvent operators $P_{\mu}^{(k)}$. As in the proof of
theorem~\ref{theointernon}, after some tedious but elementary calculations,
we\ obtain%
\[%
\begin{array}
[c]{l}%
\Vert C_{p}^{(k)}(f,g)-C^{(k)}(f,g)\Vert\\
\\
\leq c(k)~\left\{  \displaystyle\int\left\vert \left[  \eta_{p-1}^{(k-1)}%
-\pi^{(k-1)}\right]  (h)\right\vert \Upsilon_{\pi^{\left(  k\right)  }%
,\pi^{(k-1)}}^{(k),1}((f_{k},g_{k}),dh)\right. \\
\\
\hskip3cm+\left.  \displaystyle\int\left\vert \left[  \eta_{p}^{(k-1)}%
-\pi^{(k-1)}\right]  (h)\right\vert \Upsilon_{\pi^{\left(  k\right)  }%
,\pi^{(k-1)}}^{(k),2}((f_{k},g_{k}),dh)\right\}
\end{array}
\]
where $\Upsilon_{\pi^{\left(  k\right)  },\pi^{(k-1)}}^{(k),i}$, $i=1,2$, is a
pair of bounded integral operators from $\mathcal{B}(S^{(k)})^{2}$ into
$\mathcal{B}(S^{(k-1)})$ such that
\[
\int~\Vert h\Vert~\Upsilon_{\pi^{\left(  k\right)  },\pi^{(k-1)}}%
^{(k),i}((f_{k},g_{k}),dh)\leq c(k)~\Vert f_{k}\Vert\Vert g_{k}\Vert
\]
Combining the generalized Minkowski integral inequality with (\ref{dlp}) we
prove the following almost sure convergence result
\[
\lim_{p\rightarrow\infty}\Vert C_{p}^{(k)}(f,g)-C^{(k)}(f,g)\Vert=0
\]
Arguing as in the proof of theorem~\ref{theointernon}, we obtain the following
almost sure convergence result
\begin{align*}
\lim_{n\rightarrow\infty}\sum_{p=0}^{n}w_{n}^{(k)}(p)^{2}~C_{p}^{(k)}%
(f,g)(X_{p}^{(k)}) &  =\lim_{n\rightarrow\infty}\sum_{p=0}^{n}w_{n}%
^{(k)}(p)^{2}~C^{(k)}(f,g)(X_{p}^{(k)})\\
&  =\lim_{n\rightarrow\infty}\sum_{p=0}^{n}w_{n}^{(k)}(p)^{2}~\Phi^{\left(
k\right)  }\left(  \eta_{p-1}^{(k-1)}\right)  (C^{(k)}(f,g))\\
&  =\Phi^{\left(  k\right)  }\left(  \pi^{(k-1)}\right)  (C^{(k)}%
(f,g))=\pi^{(k)}(C^{(k)}(f,g))
\end{align*}
This yields the almost sure convergence
\[
\lim_{n\rightarrow\infty}\langle\mathcal{M}^{(n)}(f),\mathcal{M}%
^{(n)}(g)\rangle_{l}=\mathcal{C}_{l}^{(k)}(f,g):=\sum_{0\leq k\leq l}\pi
^{(k)}(C^{(k)}(f,g))
\]
Using the same arguments as the ones used in the proof of theorem 4.4
in~\cite{bdd08}, we conclude that the $\mathbb{R}^{d}$-valued martingale
$\mathcal{M}_{l}^{(n)}(f)$ converges in law as $n$ tends to infinity to a
martingale $\mathcal{M}_{l}(f)$ with a predictable bracket given for any air
of indexes $1\leq j,j^{\prime}\leq d$ by
\[
\langle\mathcal{M}(f^{j}),\mathcal{M}(f^{j^{\prime}})\rangle_{l}%
=\mathcal{C}_{l}^{(k)}(f^{j},f^{j^{\prime}})
\]
This ends the proof of the proposition. \hfill\hbox{\vrule height 5pt width
5pt depth 0pt}\medskip\newline

\section{Path space i-MCMC models\label{sec:pathspacemodels}}

\label{pathmodels}

The aim of this final section is to show that the multivariate CLT discussed
in section \ref{ssr} can be generalized directly to analyze the fluctuations
of the occupation measures of $(X_{n}^{(k)})_{0\leq k\leq l}$ around the
limiting tensor product measure $\otimes_{0\leq k\leq l}\pi^{(k)}$, for any
time horizon $l\geq0$. To do this, we check here that the regularity
conditions discussed in Section \ref{sec:regularityconditions} remain valid on
path spaces.

Let us fix a final time horizon $l$ for (\ref{phi}) and consider the path
space model given by
\[
\forall n\geq0\qquad X_{n}^{[l]}:=(X_{n}^{(0)},\ldots,X_{n}^{(l)})\in
S^{[l]}:=\prod_{0\leq k\leq l}S^{(k)}%
\]
For every $0\leq k\leq l$ we denote by $\eta^{(k)}\in\mathcal{P}(S^{(k)})$ the
image measure of a measure $\eta\in\mathcal{P}(S^{[l]})$ on the $k$-th
coordinate level set $S^{(k)}$ of the product space $S^{[l]}:=\prod_{0\leq
k\leq l}S^{(k)}$. Using this notation, it is easy to check that $X_{n}^{[l]}$
is an $S^{[l]}$-valued self-interacting Markov chain with transitions defined
by
\begin{equation}
\mathbb{P}(X_{n+1}^{[l]}\in dx~|~(X_{p}^{[l-1]})_{0\leq p\leq n},X_{n}%
^{[l]})=M_{\eta_{n}^{[l-1]}}^{[l]}(X_{n}^{[l]},dx)\label{pathocc}%
\end{equation}
with the occupation measures $\eta_{n}^{[l-1]}$ and the collection of
transitions $M_{\eta_{n}^{[l-1]}}^{[l]}$ defined by the following formulae
\[
\eta_{n}^{[l-1]}:=\frac{1}{n+1}\sum_{p=0}^{n}\delta_{X_{p}^{[l-1]}}%
\quad\mbox{\rm and}\quad M_{\eta_{n}^{[l-1]}}^{[l]}(X_{n}^{[l]},dx)=\prod
_{0\leq k\leq l}M_{\eta_{n}^{(k-1)}}^{(k)}(X_{n}^{(k)},dx^{l})
\]
where $x:=(x^{0},\ldots,x^{l})\in S^{[l]}$, $dx:=dx^{0}\times\cdots\times
dx^{l}$ and we have used here the convention $M_{\eta_{n}^{(-1)}}%
^{(0)}=M^{(0)}$. It is straightforward to check that (\ref{pathocc}) coincides
with the i-MCMC algorithm associated to the limiting evolution equation
\[
\pi^{\lbrack l]}=\Phi^{\left[  l\right]  }(\pi^{\lbrack l-1]})\quad\mbox{\rm
with}\quad\pi^{\lbrack l]}:=\pi^{(0)}\otimes\ldots\otimes\pi^{(l)}%
\]
and the invariant measure mapping
\[
\Phi^{\left[  l\right]  }~:~\mu\in\mathcal{P}(S^{[l-1]})\mapsto\Phi^{\left[
l\right]  }(\mu):=\pi^{(0)}\otimes\Phi^{\left(  1\right)  }(\mu^{(0)}%
)\otimes\ldots\otimes\Phi^{\left(  l\right)  }(\mu^{(l-1)})\in\mathcal{P}%
(S^{[l]})
\]

To describe the main result of this section, we need to introduce some
additional notation. For any $0\leq k_{1}\leq k_{2}$, we set
\[
S^{[k_{1},k_{2}]}:=\prod_{k_{1}\leq k\leq k_{2}}S^{(k)}\quad
\mbox{\rm  and}\quad\pi^{\lbrack k_{1},k_{2}]}:=\otimes_{k_{1}\leq k\leq
k_{2}}\pi^{(k)}\in\mathcal{P}(S^{[k_{1},k_{2}]})
\]
For any $0\leq k<l$, any pair $(\mu_{1},\mu_{2})\in\mathcal{P}(S^{[0,k]}%
)\times\mathcal{P}(S^{[k+2,l+1]})$ and any integral operator $D$ from
$S^{(k)}$ into $S^{(k+1)}$, we denote by $\mu_{1}\otimes D\otimes\mu_{2}$ the
operator from $S^{[l]}$ into $S^{[l+1]}$
\[
(\mu_{1}\otimes D\otimes\mu_{2})((x_{1},x_{2},x_{3}),dy_{1}\times dy_{2}\times
dy_{3}))=\mu_{1}(dy_{1})~D(x_{1},dy_{2})~\mu_{2}(dy_{3})
\]
where $(x_{1},x_{2},x_{3})\in S^{[l]}=\left(  S^{[0,k-1]}\otimes
S^{(k)}\otimes S^{[k+1,l]}\right)  $ and $(y_{1},y_{2},y_{3})\in
S^{[l+1]}=\left(  S^{[0,k]}\otimes S^{(k+1)}\otimes S^{[k+2,l+1]}\right)  $.

\begin{prop}
\label{keyprop}For any $l\geq0$, the mappings $\Phi^{\left[  l\right]  }$ and
the collection of Markov transition kernels $M_{\mu^{\left[  l-1\right]  }%
}^{[l]}$ satisfy the set of regularity conditions (\ref{firsto}),
(\ref{lipcna}), and (\ref{lipcnb}) as long as the corresponding conditions are
met for the marginal mappings $\Phi^{\left(  k\right)  }$ and the transitions
$M_{\mu^{\left(  k-1\right)  }}^{(k)}$ where $1\leq k\leq l$. In addition, the
mappings $\Phi^{\left[  l+1\right]  }$ satisfy the first order decomposition
(\ref{firsto}) with bounded integral operators $D_{[l+1]}$ from $S^{[l]}$ into
$S^{[l+1]}$ given by
\[
D_{[l+1]}=\pi^{\lbrack0,l]}\otimes D_{l+1}+\sum_{0\leq k<l}\pi^{\lbrack
0,k]}\otimes D_{k+1}\otimes\pi^{\lbrack k+2,l+1]}%
\]

\end{prop}

Before presenting the proof of this proposition, we emphasize that this latter
directly implies that the multivariate CLT stated in section~\ref{ftcl} is
also valid for the path space i-MCMC algorithm discussed above. In other
words, for every $k\geq0$, the sequence of random fields $(U_{n}^{[k]}%
)_{n\geq0}$ on $\mathcal{B}(S^{[k]})$ defined by
\[
U_{n}^{[k]}:=\sqrt{n}~\left[  \eta_{n}^{[k]}-\pi^{\lbrack k]}\right]
=\frac{1}{\sqrt{n+1}}\sum_{p=0}^{n}\left[  \delta_{(X_{p}^{(0},\ldots
,X_{p}^{(k)})}-\left(  \pi^{(0)}\otimes\ldots\otimes\pi^{(k)}\right)  \right]
\]
converges in law, as $n$ tends to infinity and in the sense of finite
dimensional distributions, to a sequence of Gaussian random fields $U^{[k]}$
defined as $U^{(k)}$ by replacing the semigroups $D_{l_{1},l_{2}}$ and the
limiting measures $\pi^{(l)}$ by the corresponding objects on path spaces.

\textbf{Proof of proposition~\ref{keyprop} : }\newline With some obvious
notation, we have
\[
\forall n\geq1\qquad\left(  M_{\mu^{\left[  l-1\right]  }}^{[l]}\right)
^{n}=\otimes_{0\leq k\leq l}\left(  M_{\mu^{(k-1)}}^{(k)}\right)  ^{n}%
\]
Using the fact that
\[
\Vert\otimes_{0\leq k\leq l}\mu^{\left(  k\right)  }-\otimes_{0\leq k\leq
l}\eta^{\left(  k\right)  }\Vert\leq\sum_{0\leq k\leq l}\Vert\mu^{\left(
k\right)  }-\eta^{\left(  k\right)  }\Vert
\]
for any sequence of probability measures $\mu^{\left(  k\right)  }%
,\eta^{\left(  k\right)  }\in\mathcal{P}(S^{(k)})$, with $0\leq k\leq l$, we
prove that
\[
\beta\left(  \left(  M_{\mu^{\left[  l-1\right]  }}^{[l]}\right)  ^{n}\right)
\leq\sum_{0\leq k\leq l}\beta\left(  \left(  M_{\mu^{(k-1)}}^{(k)}\right)
^{n}\right)  \leq\frac{1}{l+1}\sum_{0\leq k\leq l}m_{k}(n_{k})<1
\]
as soon as
\[
n\geq n_{[l]}:=\left(  \vee_{0\leq k\leq l}\text{ }n_{k}\right)  \times\left(
1+\frac{\log{(l+1)}}{\wedge_{0\leq k\leq l}\log{(1/m_{k}(n_{k}))}}\right)
\]
We prove the pair of regularity conditions (\ref{firsto}) and (\ref{lipcnb})
by induction on the parameter $l$. We use the notation $D_{[l]}$,
$\Xi_{\lbrack l]}$, $\varXi_{[l]}$ $m_{[l]}$, $n_{[l]}$ and $\Gamma_{\lbrack
l],\mu}$ the corresponding objects introduced in the statement of conditions
(\ref{firsto}), (\ref{lipcna}), and (\ref{lipcnb}). The results are clearly
true for $m=0$ with
\[
\Phi^{\left[  0\right]  }(\mu):=\pi^{(0)}\quad\mbox{\rm and}\quad M_{\eta
_{n}^{[-1]}}^{[0]}:=M^{(0)}%
\]
In this case, we readily find that
\[
D_{[0]}=D_{0}=0,~\Xi_{\lbrack l]}=0,~m_{[0]}=m_{0},~n_{[0]}=n_{0}%
\quad\mbox{\rm and}\quad\Gamma_{\lbrack0],\mu}=\Gamma_{0,\mu}=0
\]
Assume now that the result has been proved at some rank $l$. For any measure
$\mu$ on $S^{[l]}=S^{[l-1]}\times S^{(l)}$, we denote by $\mu^{\lbrack l-1]}$
and $\mu^{(l)}$ its image measures on $S^{[l-1]}$ and $S^{(l)}$. In this
notation we have
\[
\Phi^{\left[  l+1\right]  }(\mu)=\Phi^{\left[  l\right]  }(\mu^{\lbrack
l-1]})\otimes\Phi^{\left(  l+1\right)  }(\mu^{(l)})
\]
and
\[
M_{\mu^{\left[  l\right]  }}^{[l+1]}((u,v),dx\times dy))=M_{\mu^{\lbrack
l-1]}}^{[l]}(u,dx)\times M_{\mu^{(l)}}^{(l+1)}(v,dy)
\]
for any $(u,v)\in S^{[l+1]}=(S^{[l]}\times S^{(l)})$ and $(x,y)\in
S^{[l+1]}=(S^{[l]}\times S^{(l+1)})$. After some elementary computations,
using the decomposition
\[%
\begin{array}
[c]{l}%
\left[  \Phi^{\left[  l+1\right]  }(\mu)-\Phi^{\left[  l+1\right]  }%
(\pi^{\lbrack l]})\right] \\
\\
=\Phi^{\left[  l\right]  }(\mu^{\lbrack l-1]})\otimes\Phi^{\left(  l+1\right)
}(\mu^{(l)})-\Phi^{\left[  l\right]  }(\pi^{\lbrack l-1]})\otimes\Phi^{\left(
l+1\right)  }(\pi^{(l)})\\
\\
=\Phi^{\left[  l\right]  }(\pi^{\lbrack l-1]})\otimes\left[  \Phi^{\left(
l+1\right)  }(\mu^{(l)})-\Phi^{\left(  l+1\right)  }(\pi^{(l)})\right] \\
\\
\qquad+\left[  \Phi^{\left[  l\right]  }(\mu^{\lbrack l-1]})-\Phi^{\left[
l\right]  }(\pi^{\lbrack l-1]})\right]  \otimes\Phi^{\left(  l+1\right)  }%
(\pi^{(l)})\\
\\
\qquad\qquad+\left[  \Phi^{\left[  l\right]  }(\mu^{\lbrack l-1]}%
)-\Phi^{\left[  l\right]  }(\pi^{\lbrack l-1]})\right]  \otimes\left[
\Phi^{\left(  l+1\right)  }(\mu^{(l)})-\Phi^{\left(  l+1\right)  }(\pi
^{(l)})\right]
\end{array}
\]
we find that the first order condition (\ref{firsto}) is satisfied with an
integral operator $D_{[l+1]}$ from $S^{[l]}$ into $S^{[l+1]}$ defined for any
$f\in$ by
\begin{align*}
D_{[l+1]}(u,d(x,y)) &  =\Phi^{\left[  l\right]  }(\pi^{\lbrack l-1]}%
)(dx)~D_{(l+1)}(u,dy)+D_{[l]}(u,dx)~\Phi^{\left(  l+1\right)  }(\pi
^{(l)})(dy)\\
&  =\pi^{\lbrack l]}(dx)~D_{(l+1)}(u,dy)+D_{[l]}(u,dx)~\pi^{(l+1)}(dy)
\end{align*}
Condition (\ref{lipcnb}) is proved using the same type of arguments. This ends
the proof of the proposition. \hfill\hbox{\vrule height 5pt width 5pt
depth 0pt}\medskip\newline

\end{document}